\documentclass[11pt,a4paper]{article}
\usepackage{amsmath}
\usepackage{amsthm}
\usepackage{amsfonts}
\usepackage{multirow}
\usepackage[table,xcdraw]{xcolor}
\usepackage{subcaption}
\usepackage{float}
\usepackage{algorithm}
\usepackage{algpseudocode}
\usepackage{hyperref}
\usepackage{tabularx}
\usepackage{amssymb}
\usepackage{graphicx}
\usepackage{tcolorbox}
\usepackage{booktabs}
\usepackage[table,xcdraw]{xcolor}
\usepackage{multirow}
\usepackage{rotating,tabularx}

\usepackage[left=2.5cm,right=2.5cm,top=2.5cm,bottom=2.5cm]{geometry}
\newcolumntype{L}[1]{>{\raggedright\let\newline\\\arraybackslash\hspace{0pt}}m{#1}}
\newcolumntype{C}[1]{>{\centering\let\newline\\\arraybackslash\hspace{0pt}}m{#1}}
\newcolumntype{R}[1]{>{\raggedleft\let\newline\\\arraybackslash\hspace{0pt}}m{#1}}

\newtheorem{Theorem}{Theorem}[section]
\newtheorem{Proposition}[Theorem]{Proposition}
\newtheorem{Remark}[Theorem]{Remark}
\newtheorem{Lemma}[Theorem]{Lemma}
\newtheorem{Corollary}[Theorem]{Corollary}
\newtheorem{Definition}[Theorem]{Definition}
\newtheorem{Example}[Theorem]{Example}

\usepackage{hyperref}
\hypersetup{colorlinks=true,urlcolor=blue}
\expandafter\let\expandafter\oldproof\csname\string\proof\endcsname
\let\oldendproof\endproof
\renewenvironment{proof}[1][\proofname]{
\oldproof[\ttfamily\scshape \bf #1.]
}{\oldendproof}

\def\dom{{\rm dom}\,}

\def\epi{{\rm epi\,}}

\def\ox{\overline{x}}
\def\oy{\overline{y}}
\def\oz{\overline{z}}

\def\tto{\rightrightarrows}

\def\Bar{\overline}
\def\ra{\rangle}
\def\la{\langle}

\def\epsilon{\varepsilon}
\def\ox{\bar{x}}
\def\oy{\bar{y}}
\def\oz{\bar{z}}

\def\ov{\bar{v}}

\def\gph{\mbox{\rm gph}\,}
\def\epi{\mbox{\rm epi}\,}

\def\dom{\mbox{\rm dom}\,}

\def\O{\Omega}
\def\ph{\varphi}

\def\st{\stackrel}
\def\oR{\Bar{\R}}
\def\lm{\lambda}

\def\al{\alpha}

\def \N{{\rm I\!N}}
\def \R{{\rm I\!R}}

\def\Limsup{\mathop{{\rm Lim}\,{\rm sup}}}

\def\Limsup{\mathop{{\rm Lim}\,{\rm sup}}}

\numberwithin{equation}{section}

\title{\bf Convexity Characterizations of Functions and Sets via  Variational Modulus of Convexity}
\author{Vo Thanh Phat\footnote{Department of Mathematics and Statistics,  University of North Dakota, Grand Forks, North Dakota, USA. E-mail: thanh.vo.1@und.edu.}\and  R. Tyrrell Rockafellar\footnote{University of Washington, Department of Mathematics, Box 354350, Seattle, Washington, USA. Email:
rtr@uw.edu} }
\begin{document}  
\maketitle

\noindent
{\small{\bf Abstract}. This paper investigates the relationship between variational convexity and classical convexity of functions and sets. By employing the variational modulus of convexity, we establish new characterizations of convexity for nonsmooth functions. Building on these characterizations, we further derive derivative-based criteria for convexity that apply to nonsmooth    functions. Furthermore, we develop characterizations of locally convex closed sets by exploiting the variational convexity of their indicator functions. Unlike existing results that rely on specific geometric structures of the underlying sets, our approach requires only closedness. These findings provide new connections between variational convexity and classical convexity and contribute to the development of second-order tools in variational analysis and optimization.     \\[1ex]
 }\vspace*{-0.1in}

\section{Introduction}
The study of convexity of functions and sets plays a fundamental role in optimization theory and numerical algorithms. Convexity provides a powerful framework for analyzing the global behavior of optimization problems. In particular, for a convex optimization problem, every stationary point is a global minimizer, while convexity assumptions also provide fundamental guarantees for the convergence of several important numerical methods, including proximal point methods and proximal methods of multipliers; see, e.g., \cite{r76,r76mor}. 

\medskip 
Characterizations of convex functions through second-order generalized derivatives have been extensively investigated in the literature; see, for instance, \cite{ChieuChuongYaoYen,chieuyao11,chieuhuy,clmn16,phat26}. Related characterizations of  convex sets using second-order subdifferentials have also been developed in \cite{chieuyen}. However, most existing approaches investigate convexity directly through second-order generalized differentiation techniques, often relying on specific structural properties of the functions or sets under consideration. 

\medskip 
In this paper, we develop a new perspective on the convexity of functions and sets through the notion of variational convexity, a concept introduced recently by Rockafellar in \cite{r19}. Variational convexity provides a useful framework for studying functions that are not necessarily convex in the conventional sense, while retaining certain essential convexity properties through a suitable localization of their subgradient mappings. Since its introduction, this concept has attracted considerable attention in variational analysis and optimization. First-order characterizations and applications to variational sufficiency in optimization were developed in \cite{kkmpbook,r19,r21,roc,varprox}, while further second-order characterizations were established in both finite- and infinite-dimensional settings; see \cite{g25,KKMP23-3,kmp22convex,quadcharvar,roc24} and \cite{ar27,KKMP23-2}, respectively. Moreover, variational convexity has found applications in the analysis of optimization algorithms, including proximal point methods, augmented Lagrangian methods, progressive decoupling algorithms, and Newton-type methods; see, e.g., \cite{kmpmp,lwwy26,r19sv,r21ppm,r23,r26}.

\medskip 
A natural question arising from these developments is the following: What is the precise relationship between variational convexity and classical convexity? Although variational convexity has been extensively studied as a tool for analyzing nonconvex optimization problems, its implications for the classical convexity of functions and sets have not been fully understood. To the best of our knowledge, the present paper provides a new contribution to this direction by establishing connections between these two notions through variational and second-order analysis. In particular, our results offer a new perspective on classical convexity by exploiting recent developments in variational convexity and generalized second-order differentiation.

\medskip 
The main contributions of this paper are summarized as follows.
\begin{itemize}
    \item First, we establish conditions under which the variational convexity of a nonsmooth function implies its classical convexity. Building on these results, we obtain characterizations of convexity for nonsmooth functions in terms of the variational modulus of convexity and generalized second-order differentiation. These characterizations provide a second-order variational framework for studying classical convexity. The corresponding results are presented in Sections~\ref{sec:localconvexfunction} and \ref{sec:test}.
    \item Second, we investigate the local convexity of sets through the variational convexity of their indicator functions. In particular, we establish characterizations of locally convex sets in terms of the variational modulus of convexity and generalized second-order differentiation of indicator functions. An important feature of our approach is that it requires only the closedness of the set under consideration, without imposing additional structural assumptions on its boundary or local representation. This distinguishes our results from those in \cite{chieuyen}, where particular geometric structures, such as smooth boundaries or local epigraph representations, are considered. The corresponding results are presented in Section~\ref{sec:localset}.
\end{itemize}
The   paper is organized as follows. Section~\ref{sec:prelim} recalls the preliminary concepts and results from variational analysis that are needed throughout the paper.  Sections~\ref{sec:localconvexfunction} and \ref{sec:test} are devoted to characterizing the convexity of nonsmooth functions through variational convexity and second-order generalized differentiation. In Section~\ref{sec:localset}, we establish characterizations of the local convexity of closed sets by connecting the variational convexity to their indicator functions. Finally, Section~\ref{sec:conclusion} concludes the paper and discusses several possible directions for future research. 
 
\section{Preliminaries from Variational Analysis}\label{sec:prelim}
This section recalls some notions and notation frequently used in this paper. All our considerations are given in the space $\R^n$ with the \textit{Euclidean norm} $\|\cdot\|$ and \textit{scalar/inner product} $\langle \cdot,\cdot\rangle$ defined by
$$
\|x\|:= \left(\sum_{i=1}^n x_i^2 \right)^{\frac{1}{2}} \; \text{and } \langle x,y\rangle =\sum_{i=1}^n x_i y_i \quad \text{where } x:=(x_1,\ldots, x_n), \; y:=(y_1,\ldots,y_n).
$$
Throughout the paper, $\overline{\R} := \R \cup \{\infty\}$ denotes the extended real line, and $\N := \{1,2,\ldots\}$ denotes the set of natural numbers.  The symbols $B(x,r)$ and $\mathbb{B}(x,r)$ signify  the  open ball and closed ball  centered at $x$ with radius $r > 0$, respectively. 

\medskip 
We provide the necessary background on variational analysis and generalized differentiation, drawing on the foundational works of Mordukhovich \cite{Mordukhovich06,Mor18,Mor24} and Rockafellar and Wets \cite{Rockafellar98}. 
Given a set $\Omega\subset\R^s$ with $\oz\in\O$, the (Fr\'echet) {\em regular normal cone} to $\Omega$
at $\bar{z}\in\Omega$ is defined by \begin{equation}\label{rnc}
\widehat{N}_\Omega(\bar{z}):=\left\{v\in\R^s\;\Big|\;\limsup_{z\overset{\Omega}{\rightarrow}\bar{z}}\frac{\langle
v, z-\bar{z}\rangle}{\|z-\bar{z}\|}\le 0\right\}, \end{equation} 
where the symbol $z\overset{\Omega}{\rightarrow}\bar{z}$ indicates that $z\to\bar{z}$ with $z\in\Omega$.   The
(Mordukhovich) {\em limiting normal cones} to $\Omega$ at
$\bar{z}\in\Omega$ is defined by \begin{equation}\label{lnc}
N_\Omega(\bar{z}):=\big\{v\in\R^s\;\big|\;\exists\,z_k\st{\O}{\to}\bar{z},\;v_k\to
v\;\text{ as }\;k\to\infty\;\text{ with
}\;v_k\in\widehat{N}_\Omega(z_k)\big\}. \end{equation} 
It is important to observe that the regular normal cone defined in \eqref{rnc} is always convex, whereas the limiting normal cone in \eqref{lnc} can be nonconvex, as illustrated by the graph of $|\cdot|$ at $\bar{z}=(0,0)$. Recall that the   {\em proximal normal cone} to $\Omega $ at $\bar{z}$ is defined by
$$
N_\Omega^P(\bar{z}):=
\left\{
v\in\R^s\;\big|\;
\exists\,\sigma\geq0,\ \exists\,\delta>0
\text{ such that }
\langle v,y-x\rangle
\leq \sigma\|y-x\|^2
\quad\forall y\in\Omega\cap B(x,\delta)
\right\}.
$$ 

\medskip 
Let us consider an extended-real-valued function $\varphi:\R^n\rightarrow\overline{\R}:=(-\infty,\infty]$, and its
{\em effective domain} and {\em epigraph} defined by, respectively,
\begin{equation*}
\dom\ph:=\big\{x\in\R^n\;\big|\;\ph(x)<\infty\big\}\;\mbox{ and
}\;\epi\ph:=\big\{(x,\al)\in\R^{n+1}\;\big|\;\al\ge\ph(x)\big\}.
\end{equation*}  
We always assume that $\varphi$  is proper and lower semicontinuous (l.s.c.).  The \textit{regular subdifferential} and \textit{limiting subdifferential} of $\varphi$ at $\bar{x}\in\text{dom}\varphi$  are defined, respectively by
\begin{equation}\label{FrechetSubdifferential}
\widehat \partial\varphi\left( {\bar{x} } \right): = \big\{v\in\R^n\;\big|\;(v,-1)\in
\widehat{N}_{{\rm\text{\rm epi}}\,\varphi}\big(\bar{x},\varphi(\bar{x})\big)\big\},
\end{equation}
\begin{equation}\label{lim-sub}
\partial\varphi(\ox):=\big\{v\in\R^n\;\big|\;(v,-1)\in
N_{{\rm \text{\rm epi}}\,\varphi}\big(\bar{x},\varphi(\bar{x})\big)\big\}.
\end{equation} 
Recall that the limiting subdifferential
$\partial\varphi(\ox)$ reduces to the gradient $\{\nabla\ph(\ox)\}$
if $\ph$ is ${\cal C}^1$-smooth around $\ox$ (or merely strictly
differentiable at this point), and that $\partial\ph(\ox)$ is the
subdifferential of convex analysis if $\ph$ is convex.  The {\em symmetric subdifferential}  of $\varphi$ at $\ox$ is defined by
\begin{equation}\label{symmetric}
\partial^0 \varphi(\ox) = \partial\varphi(\ox)\cup -\partial(-\varphi)(\ox). 
\end{equation}
Unlike the positively homogeneous property of regular and limiting subdifferentials, the symmetric subdifferential is {\em fully homogeneous}; that is
\begin{equation*} 
\partial^0 (\lambda\varphi)(x) = \lambda \partial^0 \varphi(x) \quad \text{for all }\; x \in \R^n, \lambda \in \R.
\end{equation*}
Another  construction used in this paper is the proximal subdifferential of $\varphi$ at $\ox$ defined by
\begin{equation}\label{proximalgradient}
\partial_p\varphi(\ox) = \left\{ {v \in {\R^n}\, |\,\mathop {\lim \inf }\limits_{x \to \overline x } \frac{{\varphi (x) - \varphi \left( {\bar{x} } \right) - \left\langle {{v},x - \bar{x}} \right\rangle }}{{\left\| {x - \bar{x} } \right\|^2}} > -\infty} \right\}.
\end{equation}
From the above definitions, we obtain the following inclusions:
\begin{equation*} 
\partial_p\varphi(\ox) \subset \widehat{\partial}\varphi(\ox) \subset \partial\varphi(\ox) \subset \partial^0\varphi(\ox).
\end{equation*}

\medskip 
Let us now review some classes of extended-real-valued functions $\varphi\colon\R^n\to\oR$ broadly used in the paper. For $s \in \R$, the function $\varphi$ is called {\em  $s$-convex} on $\Omega$   if 
\begin{equation}\label{ineqstrongfunction}
{\varphi((1-\lambda)x +\lambda y)\le(1-\lambda)\varphi(x)+\lambda \varphi(y) -\frac{s}{2}\lambda(1-\lambda)\|x-y\|^2,}
\end{equation}
for any $x,y\in\Omega,\;\lambda\in(0,1)$. 
Specifically, the function $\varphi$ is called convex when $s = 0$, {\em $s$-strongly convex} when $s > 0$, and {\em $(-s)$-hypoconvex} when $s < 0$. In the recent work \cite{varprox}, Rockafellar introduces the notion of {\em variational $s$-convexity}, which unifies the concepts of {\em variational convexity} and {\em prox-regularity} for extended-real-valued functions in finite-dimensional spaces. Earlier contributions addressing the case $s \ge 0$ appear in \cite{r19}. Further developments and applications in  finite-dimensional  settings are provided in \cite{KKMP23-3,kmp22convex,gtdquad,quadcharvar,r21,varprox,roc24,roc}, while extensions to infinite-dimensional spaces can be found in \cite{kkmpbook,KKMP23-2}. We now recall the formal definition of variational convexity. 
\begin{Definition}[\rm variational convexity {\rm \cite{varprox}}]\label{vr} \rm Given $s\in\R$, a proper l.s.c.\ function $\varphi\colon\R^n \to \overline{\R}$ is {\em variationally $s$-convex} at $\bar{x}\in\dom\ph$ for $\bar{v} \in \partial \ph(\bar{x})$ if there exist an $s$-convex function $\widehat{\ph}$ as well as convex neighborhoods $U$ of $\bar{x}$ and  $V$ of $\bar{v}$ together with $\varepsilon > 0$ such that $\widehat{\ph} \le \ph$ on $U$ and 
\begin{equation}\label{VCneighbor}
(U \times V) \cap \gph \partial \widehat{\ph} = (U_\varepsilon \times V) \cap \gph \partial \ph, \quad\widehat{\ph}(x) = \ph(x) \text{ at common elements}\;(x,v), 
\end{equation}
where $U_\varepsilon := \{x \in U \mid \ph(x) < \ph(\bar{x}) + \varepsilon\}$. When $s \le 0$, the property is called {\em $(-s)$-level prox-regularity} of $\varphi$ at $\bar{x}$ for $\bar{v}$; we simply say that $\varphi$ is {\em prox-regular} at $\bar{x}$ for $\bar{v}$ if this holds for some $s \le 0$. 
The function $\varphi$ is {\em variationally convex} at $\bar{x}$ for $\bar{v}$ when it is $0$-level prox-regular. 
When $s>0$, $\varphi$ is {\em $s$-strongly variationally convex} at $\bar{x}$ for $\bar{v}$, and is said to be {\em strongly variationally convex} if this holds for some $s>0$.  It is {\em variationally $s$-convex} at $\ox$ (without mention of a particular $\ov$) if this
property holds for every vector $\ov \in \partial\varphi(\ox).$ 
\end{Definition} 
We can replace the set $U_\epsilon$ in Definition \ref{vr} with a standard neighborhood $U$ provided that $\varphi$ is {\em subdifferentially continuous} at $\ox$ for $\ov \in \partial\varphi(\ox)$; that is,
\begin{equation}\label{subcont}
(x^k,v^k)\in \gph \partial\varphi,\ (x^k,v^k)\to (\ox,\ov)\ \Longrightarrow\ \varphi(x^k)\to \varphi(\ox).
\end{equation}
If this property holds for every $\ov\in\partial\varphi(\ox)$, then
$\varphi$ is said to be {\em subdifferentially continuous} at $\ox$.

\medskip 
The next proposition   lists the characterizations of variationally convex functions via some properties of their limiting subdifferentials.\vspace*{-0.07in}

\begin{Proposition}[subgradient characterizations of variational convexity {\rm from \cite[Theorem 1]{varprox}}]\label{subgradientvar} Let $\varphi:\R^n\to\overline{\R}$  be l.s.c.\ with $\bar{x}\in\dom\varphi$ and $\bar{v}\in \widehat{\partial}\varphi(\bar{x})$. Then for any $s \in \R$, the following are equivalent:
\begin{itemize}
    \item[\bf (i)] $\varphi$ is variationally $s$-convex at $\ox$ for $\ov$.
    \item[\bf (ii)]  There are neighborhoods $U$ of $\bar{x}$, $V$ of $\bar{v}$, and a number $\epsilon>0$ such that  
$$ 
\langle v_1 - v_2, u_1 - u_2\rangle \geq s\|u_1-u_2\|^2 \quad \text{for all }\; (u_1,v_1),(u_2,v_2) \in \gph\partial\varphi\cap (U_\epsilon \times V),
$$
where $U_\epsilon:=\{x\in U\;|\;\varphi(x)<\varphi(\ox)+\epsilon\}$. 
\item[\bf (iii)]   There are neighborhoods $U$ of $\bar{x}$, $V$ of $\bar{v}$, and a number $\epsilon>0$ such that 
$$ 
\varphi(x)\ge\varphi(u)+\langle v, x-u\rangle +\frac{s}{2}\|x-u\|^2\quad\text{for all }\;x\in U,\;(u,v)\in\gph\partial\varphi\cap (U_\epsilon\times V). 
$$  
\end{itemize}
\end{Proposition}

To measure variational convexity, the {\em modulus of variational $s$-convexity} $\operatorname{cnv} \varphi(\bar{x}|\bar{v})$ for a proper, l.s.c. function $\varphi\colon \R^n \to \overline{\R}$ at $\bar{x}\in\dom \varphi$ with $\bar{v} \in \partial \varphi(\bar{x})$ was recently introduced in  \cite{roc24}. A negative value of $\operatorname{cnv} \varphi(\bar{x}|\bar{v})$ corresponds to prox-regularity, whereas a positive value characterizes strong variational convexity. Various explicit formulas and characterizations of this modulus were established in \cite{roc24}.
\begin{Definition}[variational modulus of  convexity  {\rm \cite{roc24}}]\label{def:modu} \rm {\em The variational modulus of  convexity} of $\varphi:\R^n \to \oR$ at $\ox$ for $\ov \in \partial\varphi(\ox)$ is defined by 
$$
\begin{aligned}
\operatorname{cnv}\varphi(\ox|\ov):=\sup\Big\{s\in\R\ \Big|\ &
\exists\,\epsilon>0, \text{ a convex neighborhood }U\times V\ni(\bar{x},\bar{v})
\text{ and }\\
& \text{an $s$-convex function } \widehat{\varphi} \text{ such that }  \eqref{VCneighbor} \;   \text{holds} \Big\}  
\end{aligned}
$$
\end{Definition}
\begin{Remark}\label{rem:cnv} \rm Let $\varphi:\R^n\to\overline{\R}$ be proper and l.s.c.,
and let $v\in\partial\varphi(x)$. By the definition of the variational
modulus of convexity in Definition~\ref{def:modu}, if $
\operatorname{cnv}\varphi(x| v)\ge s,$ 
then $\varphi$ is variationally $t$-convex at $x$ for $v$ for every
$t<s$.  This conclusion does not, in general, imply variational
convexity at level $s$, since the  supremum in the definition
of $\operatorname{cnv}\varphi(x| v)$ need not be attained. Indeed, consider the  function 
$\varphi(x)=x^3$, 
and let $\bar{x}=0$ and $\bar{v}=\varphi'(\bar{x})=0$. It is clear that  $\varphi$ is not variationally
$0$-convex at $\bar{x}$ for $\bar{v}$. Moreover, $
\operatorname{cnv}\varphi(\ox| \ov)=0.$ 
This example shows that a nonnegative  variational modulus of
convexity at a point does not necessarily imply variational convexity at
level $0$.
\end{Remark}
\medskip 
Given a proper, l.s.c. function $\varphi:\R^n\to\overline{\R}$ and a parameter value $\lambda>0$. Recall that   the {\em Moreau envelope} $e_\lambda\varphi$ and {\em proximal mapping} $\textit{\rm Prox}_{\lambda\varphi}$ are defined by, respectively,
\begin{equation}\label{Moreau}
e_\lambda\varphi(x):=\inf\Big\{\varphi(y)+\frac{1}{2\lambda}\|y-x\|^2\;\Big|\;y\in \R^n\Big\},\quad x\in\R^n,
\end{equation}	
\begin{equation}\label{ProxMapping} 
P_{\lambda}\varphi(x):= \operatorname{argmin}\Big\{\varphi(y)+ \frac{1}{2\lambda}\|y-x\|^2\;\Big|\;y\in\R^n\Big\},\quad x\in\R^n.
\end{equation}	
Moreau envelopes and proximal mappings have emerged as  effective  tools in both theoretical analysis and numerical methods for nonsmooth and nonconvex optimization, providing smooth surrogates that facilitate the study of  optimality conditions, while also enabling the design and analysis of efficient numerical algorithms; see, e.g., \cite{kkmpvn,BorisKhanhPhat,kmpmp,kmptjogo,kmptmp,kmptjogo25,phat26}.  A function $\varphi$ is said to be {\em prox-bounded} if there exists $\lambda>0$ such that $e_\lambda\varphi(x)>-\infty$ for some $x\in\R^n$. The supremum over all such $\lambda$ is the \textit{threshold} $\lambda_\varphi$ of prox-boundedness for $\varphi$.

\section{Convexity of    Functions via Variational Modulus of Convexity}\label{sec:localconvexfunction}
In this section, we establish the relationship between classical convexity and  variational convexity of nonsmooth  functions. To develop this result, we first derive several auxiliary statements. The first of these provides a characterization of convex functions defined on an open, convex set.

\begin{Proposition}\label{lm:convexfunction}  Let $O \subset \R^n$ be an open, convex set, and $\varphi:\R^n\to \oR$ be  locally Lipschitz on $O$. Then the following assertions are equivalent:
\begin{itemize}
    \item[\bf (i)] $\varphi$ is convex on $O$.
    \item[\bf (ii)] The following condition holds:
    \begin{equation}\label{eq:Frechetconvex}
    \varphi(y) \geq \varphi(x) + \langle v, y-x\rangle \quad \text{for all }\; x, y \in O, \; v \in \widehat{\partial}\varphi(x).
    \end{equation}
    \item[\bf (iii)] The following condition holds:
    \begin{equation}\label{eq:Morconvex}
    \varphi(y) \geq \varphi(x) + \langle v, y-x\rangle \quad \text{for all }\; x, y \in O, \; v \in {\partial}\varphi(x).
    \end{equation}
   \item[\bf (iv)] The following condition holds:
    \begin{equation}\label{eq:symconvex}
    \varphi(y) \geq \varphi(x) + \langle v, y-x\rangle \quad \text{for all }\; x, y \in O, \; v \in {\partial}^0\varphi(x).
    \end{equation} 
    \item[\bf (v)] $\partial^0 \varphi$ is monotone on $O$ in the sense that \begin{equation}\label{monotoneO}
\langle v_1 -v_2, x_1-x_2\rangle \geq 0 \quad \text{for all }\; x_1,x_2 \in O, v_1 \in \partial^0\varphi(x_1), v_2\in\partial^0\varphi(x_2). 
 \end{equation}
\end{itemize}
\end{Proposition}
\begin{proof}  We first verify [{\bf (i)}$\Longrightarrow${\bf (ii)}]. Indeed, we take arbitrary $x, y \in O$. Since $\varphi$ is convex on $O$, for each $t\in(0,1)$ we have
\begin{align*}
\varphi(x + t (y - x)) - \varphi(x) &= \varphi(t y + (1 - t)x) - \varphi(x) \\
&\le t \varphi(y) + (1 - t)\varphi(x) - \varphi(x)\\
&= t \varphi(y) - t \varphi(x).
\end{align*}
Hence,
$$
\frac{{\varphi(x + t (y - x)) - \varphi(x)}}{t } \le \varphi(y) - \varphi(x), \quad \forall t\in (0,1).
$$
Therefore,
\begin{equation}\label{eq:liminf}
 \mathop {\liminf }\limits_{t  \to {0^ + }} \frac{{\varphi(x + t (y - x)) - \varphi(x)}}{t } \le \varphi(y) - \varphi(x).
 \end{equation}
On the other hand, it follows from the definition of Fr\'echet subdifferentials that 
\begin{equation}\label{eq:liminf2}
 \mathop {\liminf }\limits_{t  \to {0^ + }} \frac{{\varphi(x + t (y - x)) - \varphi(x)}}{t } \geq \langle v, y-x\rangle \quad \text{for all }\; v \in \widehat{\partial} \varphi(x)
\end{equation}
Combining \eqref{eq:liminf} and \eqref{eq:liminf2}, we obtain \eqref{eq:Frechetconvex}. We next verify the implication [{\bf (ii)}$\Longrightarrow${\bf (iii)}]. To proceed, we pick any $x, y \in O$ and $v \in \partial\varphi(x)$. Then there exist sequences $x_k \to x$, $v_k \to v$ such that  $v_k \in \widehat{\partial}\varphi(x_k)$ for all $k \in \N$. Since $O$ is open, we deduce that $x_k \in O$ for all sufficiently large $k \in \N$. It follows from \eqref{eq:Frechetconvex} that 
$$
\varphi(y) \geq \varphi(x_k) + \langle v_k, y -x_k\rangle \quad \text{for all sufficiently large }\; k \in \N.
$$
Since $\varphi$ is continuous on $O$, we can let $k \to \infty$ to  obtain \eqref{eq:Morconvex}. We next prove that {\bf (iii)} implies {\bf (iv)}. To proceed, we aim to show that 
\begin{equation}\label{eq:upperFrechet}
-\widehat{\partial}(-\varphi)(x) \subset \widehat{\partial}\varphi(x)  \quad \text{for all }\; x \in O. 
\end{equation}
Fix $x\in O$. For any $v\in\partial\varphi(x)$, it follows from
\eqref{eq:Morconvex} that $v\in\widehat{\partial}\varphi(x)$.
Consequently, $
\widehat{\partial}\varphi(x)=\partial\varphi(x)$. 
Since $\varphi$ is locally Lipschitz on $O$, we have
$\partial\varphi(x)\neq\varnothing$ by
\cite[Corollary 2.25]{Mordukhovich06}, and hence
$\widehat{\partial}\varphi(x)\neq\varnothing$.
If $-\widehat{\partial}(-\varphi)(x)\neq\varnothing$, then
\cite[Proposition 1.87]{Mordukhovich06} implies that
$\varphi$ is Fr\'echet differentiable at $x$, and therefore
\eqref{eq:upperFrechet} holds. Otherwise, if
$-\widehat{\partial}(-\varphi)(x)=\varnothing$,
\eqref{eq:upperFrechet} holds trivially. We next show that 
\begin{equation}\label{eq:upperlimitM}
-{\partial}(-\varphi)(x) \subset {\partial}\varphi(x)  \quad \text{for all }\; x \in O. 
\end{equation}
Indeed, fix any $x \in O$, and take any $v \in -\partial(-\varphi)(x)$. Then there exist sequences $x_k \to x$, $v_k \to v$ such that $v_k \in -\widehat{\partial}(-\varphi)(x_k)$. Since $x_k \to x$ and \eqref{eq:upperFrechet} holds, we obtain $v_k \in \widehat{\partial}\varphi(x_k)$ for all sufficiently large $k \in \N$. This implies that $v \in \partial\varphi(x)$, which clarifies \eqref{eq:upperlimitM}. Combining the latter with the definition of symmetric subdifferentials, we deduce that 
$$
\partial^0\varphi(x) = \partial\varphi(x) \cup -\partial(-\varphi)(x) = \partial\varphi(x) \quad \text{for all }\; x \in O,  
$$
and thus {\bf (iv)} follows immediately from {\bf (iii)}.   
The implication [{\bf (iv)}$\Longrightarrow${\bf (v)}] is immediate.  
It remains to prove that {\bf (v)} implies {\bf (i)}. Take arbitrary $x_1,x_2\in O$. Since $O$ is an open, convex set, there exists an open interval $I\supset [0,1]$ such that $
x_1+t(x_2-x_1)\in O$ for all $t \in I$. We consider the function $\theta:I\to \R$ given by
$$
\theta(t):=\varphi(x_1 + t(x_2-x_1)) \quad \text{for all }\;  t\in I.
$$
By the local Lipschitz continuity of $\varphi$ on $O$,
the function $\theta$ is locally Lipschitz on $I$ . 
For any $t_1,t_2\in(0,1)$ with $t_1>t_2$, take
arbitrary $v_i\in\partial^0\varphi(x_1+ t_i(x_2-x_1))$,
$i=1,2$. By \eqref{monotoneO}, we obtain
$$
0
\leq
\langle v_1-v_2,[x_1 + t_1(x_2-x_1)]-[x_1 + t_2(x_2-x_1)]\rangle
=(t_1-t_2)
\langle v_1-v_2,x_2-x_1\rangle,
$$
which implies that 
\begin{equation}\label{eq:inner>}
\langle v_1,x_2-x_1\rangle
\geq
\langle v_2,x_2-x_1\rangle.
\end{equation}
By the limiting subdifferential chain rule in \cite[Theorem 10.6]{Rockafellar98}, we have   
$$
\partial\theta(t)
\subset 
\left\{
\langle v,x_2-x_1\rangle:
v\in\partial\varphi(x_1+t(x_2-x_1))
\right\} \quad \text{for all } \; t \in (0,1),
$$
$$
\partial(-\theta)(t)
\subset 
\left\{
\langle v,x_2-x_1\rangle:
v\in\partial(-\varphi)(x_1+t(x_2-x_1))
\right\} \quad \text{for all } \; t \in (0,1),
$$
which implies that
\begin{equation}\label{eq:chainrule1var}
\partial^0\theta(t)
\subset 
\left\{
\langle v,x_2-x_1\rangle:
v\in\partial^0\varphi(x_1+t(x_2-x_1))
\right\} \quad \text{for all } \; t \in (0,1).
\end{equation}
Combining \eqref{eq:inner>} and \eqref{eq:chainrule1var}, we deduce that 
\begin{equation}\label{monotoneoftheta}
z_1 \geq z_2 \quad \text{whenever } \; z_1 \in \partial^0\theta(t_1), \; z_2 \in\partial^0\theta(t_2),\;  0<t_2<t_1<1.
\end{equation}
Pick any $0\leq a<s<t\leq 1$. It follows from the symmetric subdifferential mean value theorem   in \cite[Theorem~4.11]{Mor18} that there exist $r\in (a,s)$ and $r^\prime \in (s,t)$ such that 
\begin{equation}\label{eq:meanvalue}
\frac{\theta(s)-\theta(a)}{s-a} = z \quad \text{and }\; \frac{\theta(t)-\theta(s)}{t-s} = z' \quad \text{for some }\; z \in \partial^0\theta(r),\; z' \in\partial^0\theta(r')
\end{equation}
Combining \eqref{monotoneoftheta} and \eqref{eq:meanvalue}, it follows that 
$$
\frac{\theta(s)-\theta(a)}{s-a}
\leq \frac{\theta(t)-\theta(s)}{t-s}.
$$
It follows from the secant-slope
characterization of convexity that $\theta$
is convex on $[0,1]$. Since $x_1, x_2$ was chosen arbitrarily on $O$, we deduce that  $\varphi$
is convex on $O$. 
    
\end{proof}
\begin{Remark}\rm Characterizations of convexity in terms of subgradient mappings have been studied extensively. In particular, Clarke's book \cite{clarke} establishes that a locally Lipschitz function on $\R^n$ is convex if and only if its Clarke subdifferential is monotone. Poliquin \cite{poli90} subsequently extended this characterization to proper, extended-real-valued functions on $\R^n$ by using proximal subgradients. Proposition~\ref{lm:convexfunction} further extends this line of results to functions defined on an arbitrary open convex subset of $\R^n$, using regular, limiting, and symmetric subdifferentials. This proposition plays an important role in establishing the main results of the subsequent sections.
\end{Remark}

The following lemma addresses the relationship between local convexity and global convexity for a function defined on a convex set.
\begin{Lemma}[\text{\rm \cite[Theorem 2.8]{ly10}}]\label{localglobal} Let $\Omega$ be a convex subset in $\R^n$, and a function $\varphi:\Omega \to \R $. Then $\varphi$ is  convex if and only if for any $z \in \Omega,$ there is $r>0$ such that $\varphi$ is convex on $B(z,r)\cap \Omega.$
\end{Lemma}
Next, recall that a set-valued mapping $S:\R^n \rightrightarrows \R^m$ is {\em outer semicontinuous} at $\ox \in \R^n$ if 
$$
\Limsup_{x\to \ox} S(x) \subset S(\ox),
$$
where the  {\em Painlev\'e-Kuratowski  limit} $\displaystyle\Limsup_{x\to \ox} S(x)$   is given by
\begin{equation*} 
\underset{x\rightarrow\bar{x}}{\Limsup}\;S(x):
=\big\{ {y\in \R ^m}\;\big|\,\exists \mbox{ sequences } x_k\to \bar x,\ y_k\rightarrow y\;\mbox {with }\;y_k\in S(x_k),\;k\in\N\big
\}.
\end{equation*}
The set-valued mapping $S$ is called {\em locally bounded} at $\ox \in \R^n$ if for some neighborhood $V$ of $\ox$, the set $S(V):=\bigcup_{x\in V}S(x)$ is bounded on $\R^m$.  
Now we are ready to derive our main results. 
\begin{Theorem}[local convexity and variational convexity of functions]\label{tm:equivalent} Let $\varphi:\R^n \to \oR$ be a proper and l.s.c. function, and  $\ox \in \dom\varphi$. Suppose that $\varphi$ is locally Lipschitz around $\ox$. Then $\varphi$ is locally convex around $\ox$ if and only if $\varphi$ is variationally convex at $\ox$.  
\end{Theorem}
\begin{proof} The implication from local convexity of $\varphi$ around $\bar{x}$ to
its variational convexity at $\bar{x}$ is immediate. Thus, it remains
only to establish the converse implication.  Due to the Lipschitz continuity of $\varphi$ around $\ox$, we deduce from \cite[Theorem 9.13]{Rockafellar98} that  $\partial\varphi(\ox)$ is nonempty, compact, and the mapping $x \mapsto \partial\varphi(x)$ is outer semicontinuous, locally bounded at $\ox$.   For each $\ov \in \partial\varphi(\ox),$ we have    $\varphi$ is variationally convex at $\ox$ for $\ov$. By Proposition \ref{subgradientvar}, there exists $\epsilon(\ov)>0$ such that 
$$ 
\varphi(x) \geq \varphi(u)+ \langle v, x-u\rangle  \quad \text{for all }\; x,u  \in B(\ox, \epsilon(\ov)), v \in \partial\varphi(u)\cap B(\ov, \epsilon(\ov)).
$$ 
Note that the collection of open sets $\{B(\ov, \epsilon(\ov))|\; \ov \in \partial\varphi(\ox) \}$ is a cover of $\partial\varphi(\ox)$. By the compactness of $\partial\varphi(\ox),$ we can find $\ov_1, \ov_2,\ldots, \ov_m \in \partial\varphi(\ox)$ and $\epsilon({\ov_1}),\epsilon({\ov_2}),\ldots, \epsilon({\ov_m})>0$ such that 
\begin{equation}\label{coversubgradient}
\partial\varphi(\ox) \subset D:= \bigcup_{i=1}^m B(\ov_i, \epsilon(\ov_i)),
\end{equation}
and for each $i = 1,2,\ldots, m,$ we have the inequalities 
\begin{equation}\label{vcineqphi}
\varphi(x) \geq \varphi(u)+ \langle v, x-u\rangle  \quad \text{for all }\; x,u  \in B(\ox, \epsilon(\ov_i)), v \in \partial\varphi(u)\cap B(\ov_i, \epsilon(\ov_i)).
\end{equation}
By \cite[Theorem 5.19]{Rockafellar98}, it follows from the outer semicontinuity, local boundedness of $\partial\varphi$ at $\ox$  and \eqref{coversubgradient} that there is $\delta>0$ such that $\partial\varphi(B(\ox,\delta))\subset D$. 
Therefore, putting $\epsilon= \text{\rm min}\{\delta,\epsilon({\ov_1}),\epsilon({\ov_2}),\ldots, \epsilon({\ov_m}) \} $ and combining with \eqref{vcineqphi}, we have 
$$
\varphi(x) \geq \varphi(u) +\langle v, x-u\rangle \quad \text{for all }\; x,u \in B(\ox,\epsilon), v\in \partial\varphi(u).
$$
Using Proposition~\ref{lm:convexfunction}, we deduce that $\varphi$ is convex on $B(\ox,\epsilon)$. The proof is complete. 
    
\end{proof}

Note that the equivalence established in Theorem~\ref{tm:equivalent} may fail if the local Lipschitz continuity assumption is weakened to mere continuity around the reference point. This is illustrated by the following example. 
\begin{Example}\rm Consider the continuous function $\varphi:\R\to\R$ defined by
$$
\varphi(x)=\sqrt{|x|} \quad \text{for all } \; x\in\R.
$$
We first observe that $\varphi$ is not locally convex around $0$. Indeed, for
$x>0$, we have
$$
\varphi''(x)=-\frac{1}{4x^{3/2}}<0.
$$
By simple computations, we have
$$
\partial\varphi(x)=
\begin{cases}
\left\{\dfrac{1}{2\sqrt{x}}\right\}& \text{if } x>0,\\[6pt]
\R & \text{if } x=0 \\[6pt]
 \left\{-\dfrac{1}{2\sqrt{-x}}\right\}&\text{if } x<0.
\end{cases}
$$
Now fix any $\bar v\in\partial\varphi(0)=\R$. Choose $\delta>0$ sufficiently small so that $
\frac{1}{2\sqrt{\delta}}>|\bar v|+1.$ 
Then, for every $x\in(-\delta,\delta)\setminus\{0\}$, we have $
|\varphi'(x)|>|\bar v|+1.$ 
It follows that
$$
\operatorname{gph}\partial\varphi
\cap\big((-\delta,\delta)\times(\bar v-1,\bar v+1)\big)
=\{0\}\times(\bar v-1,\bar v+1).
$$
Choose
$M>|\bar v|+1$. Let $
g(x):=M|x|$ for all $x \in \R$. 
Clearly, $g$ is a finite-valued convex function,  $g(x) \leq \varphi(x)$ for all $x \in (-\delta,\delta)$, and 
$g(0)=\varphi(0)=0$. Therefore, the localized
subdifferential graphs of $\varphi$ and $g$ coincide, and
the two functions agree in value on this common graph. 
Thus, $\varphi$ is variationally convex at $(0,\bar v)$. Since $\bar v\in\R$ was arbitrary, $\varphi$ is variationally convex  at $0$, despite not being locally convex around $0$. 
This shows that the local Lipschitz continuity assumption in Theorem~\ref{mainre} cannot, in general, be weakened to mere continuity.
\end{Example}

We next obtain the following characterization of convex functions via the variational modulus of convexity.
\begin{Theorem}[convexity of locally Lipschitz functions via the variational modulus of  convexity] \label{mainre} Let $O \subset \R^n$ be an open, convex set, and $\varphi:\R^n \to \oR$ be locally Lipschitz on $O$. Then $\varphi$ is convex on $O$ if and only if   for any $x \in O$, $v \in \partial\varphi(x)$, we have $\operatorname{cnv}\varphi(x|v) \geq 0$. 
\end{Theorem}
\begin{proof} If $\varphi$ is convex on $O$, then $\varphi$ is variationally convex at any point $x \in O$, and thus $\operatorname{cnv}\varphi(x|v) \geq 0$ for all $v \in \partial\varphi(x)$.  Conversely, suppose that $\operatorname{cnv}\varphi(x|v) \geq 0$ for any $x \in O$ and $v \in \partial\varphi(x)$, we show that $\varphi$ is convex on $O$. To proceed, we fix any $s <0$, and consider the function $\varphi_s:\R^n \to \oR$ given by 
$$
\varphi_s(x) := \varphi(x)-\frac{s}{2}\|x\|^2 \quad \text{for all }\; x\in \R^n. 
$$
Pick any $\ox \in O$, we first show that $\varphi_s$ is convex on a neighborhood of $\ox$. Due to Remark~\ref{rem:cnv}, the condition $\operatorname{cnv}\varphi(\ox|\ov) \geq 0$ for all $\ov \in \partial\varphi(\ox)$ implies that $\varphi$ is variationally $s$-convex at $\ox$. It follows from \cite[Lemma 4.2]{kmp22convex} that $\varphi_s$ is variationally convex at $\ox$.  Combining the latter with the local Lipschitz continuity of $\varphi_s$ around $\ox$, we deduce from Theorem~\ref{tm:equivalent} that $\varphi_s$ is locally convex around $\ox$.  
Since $\ox$ is chosen arbitrarily on $O$, this means that $\varphi_s$ is locally convex on $O$. Furthermore,  follows from Lemma \ref{localglobal} that $\varphi_s$ is convex on $O$, which implies that
$$
\varphi(tx + (1-t)y) \leq t \varphi(x) + (1-t)\varphi(y) - \frac{s}{2}t(1-t)\|x-y\|^2 \quad \text{for all }\; x, y\in O, t \in [0,1].
$$
Since $s<0$ is chosen arbitrarily, we can let $s \to 0^-$ to obtain 
$$
\varphi(tx + (1-t)y) \leq t \varphi(x) + (1-t)\varphi(y)  \quad \text{for all }\; x, y\in O, t \in [0,1],
$$
and thus $\varphi$ is convex on $O$. The proof is complete. 
\end{proof}

The following counterexample demonstrates that, in the absence of local Lipschitz continuity of $\varphi$ in Theorem \ref{mainre}, the nonnegativity of the modulus of variational convexity does not ensure that the function is convex.
\begin{Example}\label{l0}\rm Consider the continuous function $\varphi:\R\to\R$
defined by
\[
\varphi(x):=-\sqrt{|x|}.
\]
For every $x\neq0$, we have
$$
\varphi'(x)=
\begin{cases}
\dfrac{1}{2\sqrt{-x}},&x<0,\\[5pt]
-\dfrac{1}{2\sqrt{x}},&x>0,
\end{cases}
\quad \text{and }\; 
\varphi''(x)=\frac{1}{4|x|^{3/2}}>0.
$$
Hence, $\varphi$ is locally strictly convex around every
$x\neq0$. 
On the other hand, for every $v\in\R$,
$$
\liminf_{x\to0}
\frac{\varphi(x)-\varphi(0)-vx}{|x|}
=
\liminf_{x\to0}
\left(-\frac{1}{\sqrt{|x|}}-
v\frac{x}{|x|}\right)
=-\infty.
$$
It follows that $
\widehat{\partial}\varphi(0)=\varnothing$. 
Moreover, since
$$
|\varphi'(x)|=\frac{1}{2\sqrt{|x|}}
\to+\infty
\quad\text{as }x\to0,
$$
we also have $\partial\varphi(0)=\varnothing$.
Consequently, $\varphi$ is variationally convex at every
point of $\operatorname{gph}\partial\varphi$, since its
subdifferential graph contains only points with $x\neq0$,
where $\varphi$ is locally strictly convex. In particular,
$$
\operatorname{cnv}\varphi(x| v)\geq0
\quad
\text{for all }\;  (x,v)\in\operatorname{gph}\partial\varphi.
$$
Furthermore, $\varphi$ is subdifferentially continuous
everywhere on its subdifferential graph.
Nevertheless, $\varphi$ is not convex, since
$$
\varphi(0)=0>
\frac{\varphi(-1)+\varphi(1)}{2}=-1.
$$
\end{Example}

\begin{Corollary}  \label{co:mainre} Let $O \subset \R^n$ be an open, convex set, a real number $s \in \R$, and $\varphi:\R^n \to \oR$ be locally Lipschitz on $O$. Then $\varphi$ is $s$-convex on $O$ if and only if   for any $x \in O$, $v \in \partial\varphi(x)$, we have $\operatorname{cnv}\varphi(x|v) \geq s $.  
\end{Corollary} 
\begin{proof} Consider the $s$-quadratic shift function $\varphi_s:\R^n \to \oR$ given by 
\begin{equation}\label{eq:s-shift}
\varphi_s(x) := \varphi(x) - \frac{s}{2}\|x\|^2 \quad \text{for all }\; x\in \R^n. 
\end{equation}
The elementary sum rule for limiting subdifferentials from \cite[Proposition
1.30]{Mor18} gives us
$$
\partial\varphi_s(x) = \partial\varphi(x) - sx \quad \text{for all }\; x \in \R^n.
$$
Therefore, for any $x \in O, v\in\partial\varphi_s(x)$, we have $v +sx \in \partial\varphi(x)$. It follows from \cite[Proposition~2.3]{roc24} that 
$$
\operatorname{cnv}\varphi_s(x|v) \geq \operatorname{cnv}\varphi(x|v+sx) - s  \geq 0 \quad \text{for all }\; x \in O, v\in\partial\varphi_s(x). 
$$
Using Theorem~\ref{mainre}, we deduce that $\varphi_s$ is convex on $O$, which tells us that $\varphi$ is $s$-convex on $O$. The proof is complete. 

\end{proof}

For a general extended-real-valued function $\varphi\colon\R^n\to\oR$, which is not necessarily locally Lipschitz, we establish an additional characterization of convexity in terms of the variational modulus of convexity. To proceed, we first recall the following lemma. 
\begin{Lemma}[convexity characterization   via Moreau envelopes {\rm from \cite[Theorem 3.17]{wang10}}]\label{convexMoreau}  Suppose that  $\varphi:\R^n\to\overline{\R}$ is a proper, l.s.c., and prox-bounded function  with
 threshold $\lambda_\varphi>0$. Then $\varphi$ is convex if and only if  the Moreau envelope $e_\lambda \varphi$ is convex for some $\lambda\in(0,\lambda_\varphi)$.  
\end{Lemma}
Now we are ready to derive our aforementioned result. 
\begin{Theorem}[convexity of hypoconvex functions via the variational modulus of  convexity]\label{thm:hypo}  Let $\varphi:\R^n\to\overline{\R}$ be a proper, l.s.c., hypoconvex function, and let $s \in \R$.  Then $\varphi$ is $s$-convex if and only if $\operatorname{cnv}\varphi(x|v) \geq s$ for all $x \in \R^n$, $v \in \partial\varphi(x)$. 
\end{Theorem}
\begin{proof} If $\varphi$ is $s$-convex, then it is variationally $s$-convex at every point in $\R^n$. Consequently, $
\operatorname{cnv}\varphi(x|v)\geq s$ for all $x \in \R^n$, $v \in \partial\varphi(x)$. 
To establish the converse implication, we first consider the case $s=0$. Indeed, it follows from \cite[Lemma~2.5]{dd22} that the hypoconvexity of $\varphi$ guarantees the existence of $\lambda_0>0$ such that, for every $\lambda\in(0,\lambda_0)$, all the assumptions and conclusions of \cite[Lemma~3.2]{ar26} are satisfied. By \cite[Theorem~3.4]{ar26}, for any $x\in\R^n$ and $v\in\partial\varphi(x)$, the Moreau envelope $e_\lambda\varphi$ is $\mathcal{C}^{1,1}$-smooth in a neighborhood of $x_\lambda:=x+\lambda v$. Moreover, we have
$$
\operatorname{cnv}e_\lambda\varphi\big(x_\lambda|\nabla e_\lambda\varphi(x_\lambda)\big)
=\frac{\operatorname{cnv}\varphi(x| v)}{1+\lambda \operatorname{cnv}\varphi(x| v)} \geq 0,
$$
which implies that $e_\lm\varphi$ is convex by using Theorem~\ref{mainre}. Combining the latter with Lemma~\ref{convexMoreau}, $\varphi$ is convex. If $s\ne 0$, we consider the $s$-quadratic shift function $\varphi_s$ defined in \eqref{eq:s-shift} and proceed similarly to the proof of Corollary~\ref{co:mainre}. 
\end{proof}
\begin{Remark}[on the proof of Theorem~\ref{thm:hypo}] \rm The main idea behind the proof of Theorem~\ref{thm:hypo} is to exploit the equivalence between the convexity of Moreau envelopes and that of the underlying functions, together with the Moreau-spectral recovery formula for the variational modulus of convexity established in \cite[Theorem~3.4]{ar26}. These results are available under prox-boundedness and prox-regularity assumptions. This raises a natural question: why do we assume hypoconvexity rather than the weaker condition of prox-regularity? Although the formula in \cite[Theorem~3.4]{ar26} requires only prox-regularity, the admissible parameter $\lambda$ for the Moreau envelope depends locally on the reference point and subgradient. Consequently, prox-regularity alone does not guarantee the existence of a uniform parameter $\lambda$ for which the Moreau envelope is sufficiently smooth at every point. This motivates our use of hypoconvexity, which provides a global counterpart of prox-regularity and ensures the uniform availability of a suitable parameter $\lambda$. Furthermore, even if we assume that $\varphi$ is prox-regular at every point in $\R^n$, the conclusion of Theorem~\ref{thm:hypo} may fail, as illustrated by Example~\ref{l0}.
\end{Remark}
The following example demonstrates that the global hypoconvexity assumption in Theorem~\ref{thm:hypo} cannot, in general, be replaced by local hypoconvexity, even when the local hypoconvexity modulus is uniform.

\begin{Example}\rm
Let $\varphi:\R\to\overline{\R}$ be defined by
$$
\varphi(x):=\delta_{\{-1,1\}}(x)
=
\begin{cases}
0, & \text{if } x\in\{-1,1\},\\
+\infty, & \text{otherwise}.
\end{cases}
$$
Then $\varphi$ is proper and l.s.c., but it is not convex, since
$\operatorname{dom}\varphi=\{-1,1\}$ is not a convex set. 
On the other hand, $\varphi$ is locally convex around every point in its domain. Consequently,
$$
\operatorname{cnv}\varphi(x| v)\geq 0
\quad\text{for all }(x,v)\in\operatorname{gph}\partial\varphi.
$$
Moreover, $\varphi$ is locally hypoconvex around every point in its domain, with the same modulus $s=0$. This example shows that local hypoconvexity, even with a uniform modulus, is not sufficient to replace the global hypoconvexity assumption in Theorem~\ref{thm:hypo}. 
\end{Example}

\section{Derivative Tests for Convexity of Nonsmooth   Functions}\label{sec:test}
In this section, we provide second-order derivative tests for nonsmooth convex functions. To proceed, we first recall the definitions of certain second-order generalized derivative structures. Let $\ph: \R^n \to \overline{\R}$ be a proper, l.s.c. function. For each $x \in \R^n$, and $v \in \R^n$, the \textit{second-order difference quotient} $\Delta_\tau^2 \varphi(x|v)$ for $\tau > 0$ is defined as
$$
    \Delta_\tau^2 \ph(x|v)(w) := \frac{\ph(x + \tau w) - \ph(x) - \tau \la v, w \ra}{\frac{1}{2}\tau^2}  \quad \text{for all }\; w \in \R^n.
$$
The {\em second subderivative} of $\varphi$ at $\ox$ for $\ov\in\R^n$ is defined as
\begin{equation*}
d^2\varphi(\ox|\ov)(w):=\liminf_{\tau\downarrow 0, w'\to w}\Delta^2_\tau\varphi(\ox,\ov)(w').
\end{equation*}
The \textit{strict second-order subderivative} of $\ph$ at $\bar{x} \in \dom \ph$ for $\bar{v} \in \R^n$ is defined as
\begin{equation}\label{strict2nd}
d^2_* \ph(\bar{x}|\bar{v})(w) = \liminf\limits_{\substack{t \downarrow 0, w' \to w \\ (x,v)\xrightarrow{\text{\rm gph} \partial \ph} (\bar{x},\bar{v}) \\ \ph(x) \to \ph(\bar{x})}} \Delta_t^2 \ph(x|v)(w').
\end{equation}
Now we are ready to obtain the following characterization of convex functions. 
\begin{Theorem}[characterization  via strict second subderivatives]\label{convexvia2nd} Let $\varphi:\R^n\to\overline{\R}$ be a proper, l.s.c. function, $O \subset \R^n$ be an open, convex set, and let $s \in \R$.  If $\varphi$ is $s$-convex on $O$, then the following condition holds: 
\begin{equation}\label{PSD2ndstrict}
    d^2_* \varphi(x|v)(w) \geq s\|w\|^2 \quad \text{for all }\; x\in O, v \in \partial\ph(x), w \in \R^n. 
\end{equation}
The converse implication also holds under any one of the following conditions:
\begin{itemize}
    \item[\bf (i)] $\varphi$ is locally Lipschitz on $O$.
    \item[\bf (ii)] $\varphi$ is hypoconvex on $O$, and $O = \R^n$. 
\end{itemize}
\end{Theorem}
\begin{proof} ($\Longrightarrow$): Suppose that $\varphi$ is $s$-convex on $O$, then $\varphi$ is variationally $s$-convex at any point $x \in O$. Pick any $x \in O$, $v \in \partial\varphi(x)$, it follows from  Proposition~\ref{subgradientvar} that there exist $\epsilon>0$,  neighborhoods $U$ of $x$, $V$ of $v$  such that for all $(x',v') \in \gph \partial\varphi$, $\varphi(x')<\varphi(x) + \epsilon$, and  all $y \in U$ such that 
$$
\varphi(y) \ge \varphi(x') + \la v', y - x' \ra + \frac{s}{2}\|y - x'\|^2.  
$$
Fix any $w \in \R^n$. For all $t > 0$, $\forall w' \in \R^n$ such that $\|t w'\|$ is sufficiently small, for all $(x',v') \in \gph \partial\varphi$, $\varphi(x')<\varphi(x) + \epsilon$, we then obtain
$$ 
\Delta_t^2 \varphi(x'|v')(w') \ge s\|w'\|^2. $$
Therefore,
 $$
 d_*^2 \varphi(x|v)(w) = \liminf\limits_{\substack{t \downarrow 0, w' \to w \\ (x',v')\xrightarrow{\text{\rm gph} \partial \varphi} (x,v) \\ \varphi(x') \to \varphi(x)}} \Delta_t^2 \varphi(x'|v')(w') \geq s\|w\|^2,
 $$
which verifies \eqref{PSD2ndstrict}. 

\medskip\noindent
($\Longleftarrow$): 
We now assume that \eqref{PSD2ndstrict} holds and that one of conditions {\bf (i)} or {\bf (ii)} is satisfied. We aim to show that $\varphi$ is $s$-convex on $O$. It follows from Theorem 2.2 in \cite{roc24} that the modulus of variational convexity can be expressed in terms of strict second subderivatives as follows:
\begin{equation}\label{cnv2nd}
\operatorname{cnv}\ph(x|v) =\inf_{\|w\|=1}d^2_*\ph(x|v)(w), \quad \text{where }\; x \in \R^n, v \in \partial\varphi(x). 
\end{equation}
For every $x \in O$ and $v \in \partial\varphi(x)$, it follows immediately from \eqref{cnv2nd} that $\operatorname{cnv}\varphi(x|v) \geq s$. Using Corollary~\ref{co:mainre} and Theorem~\ref{thm:hypo}, we deduce that $\varphi$ is $s$-convex on $O$. The proof is complete. 
    
\end{proof}

\medskip Based on the limiting normal cone     \eqref{lnc} to the graphs, we  can define the  {\em limiting
coderivative} of a set-valued mapping $F:\R^n \rightrightarrows \R^m$ at $(\ox,\oy)\in\gph F$ by 
\begin{equation}\label{lim-cod}
D^*F(\ox,\ov)(v):=\big\{u\in\R^n\;\big|\;(u,-v)\in N_{{\rm
gph}\,F}(\ox,\ov)\big\},\quad v\in\R^m. 
\end{equation}  
Recall that $F$ is called {\em positively homogeneous} if 
$$
0 \in F(0) \quad \text{and }\; F(\lm w) = \lm F(w) \quad \forall w\in \R^n,\lm>0. 
$$
Next, we recall the {\em second-order subdifferential} \cite{m92}, or {\em generalized Hessian}, $\partial^2\ph(\ox,\ov)\colon\R^n\tto\R^n$ of $\ph\colon\R^n\to\oR$ at $\ox\in\dom\ph$ for $\ov\in\partial\ph(\ox)$ as the coderivative of the subgradient mapping
\begin{equation}\label{2nd}
\partial^2\ph(\ox,\ov)(u):=\big(D^*\partial\ph\big)(\ox,\ov)(u)\;\mbox{ for all }\;u\in\R^n.
\end{equation}
If $\ph$ is ${\cal C}^2$-smooth around $\ox$, then we have
$$ 
\partial^2\ph(\ox)(u)=\big\{\nabla^2\ph(\ox)u\big\}\;\mbox{ for all }\;u\in\R^n.
$$
Note also that $\partial^2\varphi(\ox,\ov)$ is positively homogeneous. As in \cite{g25}, second-order subdifferentials have been modified for adaptation to work best with subdifferential
mappings. Gfrerer \cite{g25} adds the requirement that $\varphi(x_k) \to \varphi(\ox)$ in the definition of second-order subdifferential to replace \eqref{2nd} by, in his notation, $D^*_\varphi (\partial\varphi)$. More specfically, 
\begin{equation}\label{eq:modified2nd}
D^*_\varphi (\partial\varphi)(\ox,\ov)(v):= \big\{u\in\R^n\;\big|\;(u,-v)\in N^\varphi_{{\rm
gph}\,\partial\varphi}(\ox,\ov)\big\},\quad v\in\R^n,
\end{equation}
where 
$$
\begin{aligned}
N^\varphi_{\operatorname{gph}\partial\varphi}(\bar x,\bar v)
:=\big\{(u,v)|\;&\exists\,(x_k,v_k)\xrightarrow{\operatorname{gph}\partial\varphi}(\bar x,\bar v),\\
& (u_k,v_k)\to(u,v),\ \varphi(x_k)\to\varphi(\bar x),\
(u_k,v_k)\in\widehat N_{\operatorname{gph}\partial\varphi}(x_k,v_k)\big\}.
\end{aligned}
$$
Note that if $\varphi$ is subdifferentially continuous at $\bar x$ for $\bar v\in\partial\varphi(\bar x)$, then the two definitions in \eqref{2nd} and \eqref{eq:modified2nd} coincide.
Characterizations of important classes of nonsmooth convex functions, including $\mathcal{C}^{1,1}$-smooth, $\mathcal{C}^1$-smooth, piecewise $\mathcal{C}^2$-smooth, and lower-$\mathcal{C}^2$ functions, have been extensively studied in the literature \cite{ChieuChuongYaoYen,chieuyao11,chieuhuy,clmn16,phat26} using second-order subdifferentials \eqref{2nd}. However, as demonstrated in \cite{chieuhuy}, such characterizations cannot, in general, be extended to arbitrary locally Lipschitz functions. In contrast, Theorem~\ref{convexvia2nd} provides a new characterization of convexity for locally Lipschitz functions in terms of strict second subderivatives \eqref{strict2nd}. The following example illustrates the distinction between these two approaches by exhibiting a locally Lipschitz nonconvex function whose limiting second-order subdifferential satisfies the positive-semidefinite property everywhere, whereas its strict second subderivative takes the value $-\infty$. This example demonstrates that strict second subderivatives can detect nonconvexity in cases where the positive-semidefinite property of limiting second-order subdifferentials fails to do so.
\begin{Example} \rm Let $E\subset\R$ be a measurable set such that, for every
nonempty open interval $I\subset\R$, both $I\cap E$ and
$I\cap E^c$ have positive Lebesgue measure. Such sets exist
because of  \cite[p.~307]{s81}. We consider the function $\varphi:\R\to\R$ defined by
$$
\varphi(x):=\int_0^x\chi_E(t)\,dt\quad \text{for all }\; x \in \R.
$$
As established in \cite[Theorem 4.2]{chieuhuy}, $\varphi$ is Lipschitz continuous
and
\[
\partial\varphi(x)=[0,1]
\qquad\text{for every }x\in\R.
\]
Moreover, its limiting second-order subdifferential $\partial^2\varphi$ satisfies
the positive-semidefinite property in the sense that
$$
\langle z, w\rangle \geq 0 \quad \text{for all }\; z \in \partial^2\varphi(x,v)(w),\; v\in\partial\varphi(x),\; w\in \R,
$$
although $\varphi$ is not
convex. 
We next show that its strict second subderivative takes a negative value. More precisely, we claim that
$$
d_*^2\varphi\left(\bar{x}\big|\frac12\right)(1)
=-\infty
\qquad\text{for every }\bar{x}\in\R.
$$
Fix an arbitrary $\bar{x}\in\R$.
Since $E^c$ has positive measure in every nonempty open interval,
the Lebesgue density theorem guarantees the existence of a
sequence of density points $
x_k\in E^c$ for all $k \in \N$ and $x_k\to\bar{x}$ as $k \to \infty$. 
For each $k \in \N$, since $x_k$ is a density point of $E^c$, we have
$$
\lim_{t\downarrow0}
\frac{1}{t}\int_{x_k}^{x_k+t}\chi_E(s)\,ds=0.
$$
Consequently, we can choose $t_k>0$ sufficiently small such that
$$
0<t_k<\frac1k
\quad \text{and }\;
\frac{1}{t_k}\int_{x_k}^{x_k+t_k}\chi_E(s)\,ds
<\frac1k,
$$
and thus 
$t_k\downarrow 0$. 
Now choose $
v_k:=\frac12$ and $u_k:=1$ for all $k \in \N$. 
Since $\partial\varphi(x_k)=[0,1]$, we have
$$
v_k\in\partial\varphi(x_k) \; \text{for all } k \in \N, \; \lim_{k\to \infty} v_k = \frac{1}{2},\; \lim_{k\to\infty} u_k = 1, \; \lim_{k\to \infty}\varphi(x_k) = \varphi(\ox)
$$
Moreover, we have 
\begin{align*}
\frac{2\big[
\varphi(x_k+t_ku_k)-\varphi(x_k)-t_kv_ku_k
\big]}{t_k^2}
\quad & =
\frac{2\big[
\varphi(x_k+t_k)-\varphi(x_k)-\frac12t_k
\big]}{t_k^2} \\
& =
\frac{2}{t_k}
\left[
\frac{\varphi(x_k+t_k)-\varphi(x_k)}{t_k}
-\frac12
\right] \\
& =
\frac{2}{t_k}
\left[
\frac{1}{t_k}\int_{x_k}^{x_k+t_k}
\chi_E(s)\,ds-\frac12
\right]
<
\frac{2}{t_k}\left(\frac1k-\frac12\right).
\end{align*}
For every $k\ge3$, we have $
\frac1k-\frac12<-\frac16.$ 
It follows that
$$
\frac{2\big[
\varphi(x_k+t_k)-\varphi(x_k)-\frac12t_k
\big]}{t_k^2}
<-\frac{1}{3t_k}
\longrightarrow-\infty,
$$
because $t_k\downarrow0$. 
Therefore, by the definition of the strict second subderivative, we have 
$$
d_*^2\varphi\left(\bar{x}\big|\frac12\right)(1)
=-\infty.
$$
This demonstrates that the strict second subderivative can
detect nonconvexity in this example, while the positive-semidefinite
property of the limiting second-order subdifferential cannot.
\end{Example}

While the second-order subdifferential \eqref{2nd} does not generally characterize convexity for locally Lipschitz functions, it does so under the additional assumption of prox-regularity.
\begin{Theorem}[characterization   via  second-order subdifferentials]\label{thm:lowerC2}  Let $\varphi:\R^n\to\overline{\R}$ be a proper, l.s.c. function, $O \subset \R^n$ be an open, convex set, and let $s \in \R$.   If $\varphi$ is $s$-convex on $O$, then the following condition holds:
\begin{equation}\label{eq:codebased} 
\langle z, w\rangle \geq s\|w\|^2 \quad \text{for all }\; z \in D^*_\varphi (\partial\varphi)(x,v)(w),\; x \in O, \;v\in\partial\varphi(x),\;w\in\R^n.
\end{equation}
The converse implication also holds under any one of the following conditions:
\begin{itemize}
    \item[\bf (i)]  $\varphi:\R^n \to \oR$ be locally Lipschitz   prox-regular on $O$. 
    \item[\bf (ii)] $\varphi$ is hypoconvex on $O$, and $O=\R^n$.  
\end{itemize}
\end{Theorem}
\begin{proof}  It follows from \cite[Theorem~3.8]{roc24} that if $\varphi$ is prox-regular on $O$,  then the modulus of variational convexity can be expressed in terms of \eqref{eq:modified2nd} as follows:
\begin{equation}\label{eq:2ndviavariational}
\operatorname{cnv}\ph(x|v)  = \operatorname{def}  D^*_\varphi (\partial\varphi)(x,v) \quad \text{for all }\; x \in O, \; v \in \partial\varphi(x),
\end{equation}
where 
$$
\operatorname{def}  D^*_\varphi (\partial\varphi)(x,v) := \inf \{\langle z,u\rangle |\;  z \in  D^*_\varphi (\partial\varphi)(x,v)(u),\; \|u\|=1\}.
$$
We first verify the implication ``($\Longrightarrow$)''. Indeed, suppose that $\varphi$ is $s$-convex on $O$, then $\varphi$  is variationally $s$-convex at any point $x \in O$, which implies that $\operatorname{cnv}\varphi(x|v) \geq s$ for all $x \in O$ and  $v\in \partial\varphi(x)$. Pick any $x \in O$, $v \in \partial\varphi(x)$, it follows from  \eqref{eq:2ndviavariational} that 
\begin{equation}\label{2ndnorm1}
\langle z, u\rangle \geq s  \quad \text{for all }\; z \in  D^*_\varphi (\partial\varphi)(x,v)(u), \; \|u\|=1. 
\end{equation}
Pick any $w \in \R^n \setminus \{0\}$, $z \in D^*_\varphi (\partial\varphi)(x,v)(w)$, and define $u:= w/\|w\|$. It follows from the fact that $ D^*_\varphi (\partial\varphi)(x,v)$ is positively homogeneous  that $z/\|w\| \in D^*_\varphi (\partial\varphi)(x,v)(u)$.   Combining the latter with the fact that $\|u\|=1$ and \eqref{2ndnorm1}, we must have 
$$
\left\langle \frac{z}{\|w\|}, u  \right\rangle \geq s \iff \left\langle \frac{z}{\|w\|}, \frac{w}{\|w\|}  \right\rangle \geq s \iff \langle z , w\rangle \geq s\|w\|^2,
$$
which verifies \eqref{eq:codebased} in the case where $w \ne 0$. In the case where $w =0$, \eqref{eq:codebased} holds trivially. 

\medskip\noindent
($\Longleftarrow$): We now assume that \eqref{eq:codebased} holds and that one of conditions {\bf (i)} or {\bf (ii)} is satisfied. We aim to show that $\varphi$ is $s$-convex on $O$.   
For every $x \in O$ and $v \in \partial\varphi(x)$, it follows immediately from \eqref{eq:2ndviavariational} that $\operatorname{cnv}\varphi(x|v) \geq s$. Using Corollary~\ref{co:mainre} and Theorem~\ref{thm:hypo}, we deduce that $\varphi$ is $s$-convex on $O$. The proof is complete. 
    
\end{proof}
\begin{Remark}\rm Note that, for a locally Lipschitz function $\varphi$ on $O$, prox-regularity on $O$ is equivalent to being lower-$\mathcal{C}^2$ on $O$; see \cite[Proposition 13.33]{Rockafellar98}. Recall that $\varphi$ is said to be {\em lower-$\mathcal{C}^2$} on $O$ if, for every $\bar{x}\in O$, there exists a neighborhood $V$ of $\bar{x}$ such that
$$
\varphi(x)=\max_{t\in T}\varphi_t(x),
$$
where each function $\varphi_t$ is of class $\mathcal{C}^2$ on $V$, the index set $T$ is compact, and the functions $\varphi_t$, together with their first- and second-order partial derivatives, depend continuously on $(x,t)\in V\times T$. Hence, Theorem~\ref{thm:lowerC2} provides a characterization of convexity for lower-$\mathcal{C}^2$ functions on arbitrary open convex sets. This extends the results of \cite[Theorem 4.1 and Corollary 4.3]{clmn16}, which characterize convexity on the entire space $\R^n$ and consider only the case $s\geq 0$.
\end{Remark}

We next employ the {\em strict graphical derivative} of subgradient mappings  to characterize the convexity of nonsmooth functions. This derivative provides a robust primal generalized-derivative framework, in contrast to the {\em graphical derivative}, which does not possess the same robustness property as the coderivative \eqref{lim-cod} in the dual setting. To recall these notions, we recall some concepts.  
For $x \in \R^n$, $v \in \R^m$, the {\em difference quotient mapping} of $F$ at $(x,v) \in \R^n\times \R^m$ is    defined by
$$
\Delta_\tau F(x|v)(w):= \frac{F(x+\tau w)-v}{\tau} \quad \text{for }\; \tau >0, w \in \R^n.
$$
The {\em graphical derivative} and {\em strict graphical derivative} of $F$ at $(\ox,\ov)$ are defined, respectively, by 
\begin{equation}\label{graphical}
\gph DF(\ox,\ov) = \Limsup_{\tau \searrow 0 } \gph \Delta_\tau F(\ox|\ov),
\end{equation}
\begin{equation}\label{strictgraphical}
\gph D_* F(\ox,\ov) = \Limsup\limits_{\substack{t \searrow 0 \\ (x,v)\xrightarrow{\text{\rm gph F}}(\ox,\ov)}} \Delta_\tau F(x|v).
\end{equation}
To better capture the behavior of subgradient mappings, Rockafellar \cite{roc24} introduced a modification of the strict graphical derivative that incorporates information about the underlying function, which he termed the \emph{strict second-order subdifferential}. More precisely,
\begin{equation}\label{eq:2ndstrict}
\gph \partial^2_* \varphi(\ox,\ov):= \Limsup\limits_{\substack{t \searrow 0 \\ (x,v)\xrightarrow{\text{\rm gph}\partial\varphi}(\ox,\ov)\\\varphi(x)\to \varphi(\ox) }} \Delta_\tau \partial\varphi(x|v).
\end{equation}
Note that if $\varphi$ is subdifferentially continuous at $\bar x$ for $\bar v$, then the strict graphical derivative \eqref{strictgraphical} of the subgradient mapping and the strict second-order subdifferential \eqref{eq:2ndstrict} coincide; that is,
$$
\gph \partial^2_* \varphi(\ox,\ov) =\gph (D_* \partial\varphi)(\ox,\ov).
$$
Characterizations of nonsmooth convex functions in terms of the graphical derivative \eqref{graphical} have been investigated in \cite{phat26} using Moreau envelope techniques and under an additional  hypoconvexity assumption. In the present work, we employ the strict second-order subdifferential \eqref{eq:2ndstrict} to characterize the convexity of nonsmooth functions under  weaker assumptions than those required by the second-order subdifferential criteria in Theorem~\ref{thm:lowerC2}. Here is our result.
\begin{Theorem}[characterization via strict second-order subdifferentials] \label{thm:proximalsub}  Let $\varphi:\R^n\to\overline{\R}$ be a proper, l.s.c. function, $O \subset \R^n$ be an open, convex set, and let $s \in \R$.   If $\varphi$ is $s$-convex on $O$, then the following condition holds: 
\begin{equation}\label{eq:graphbased} 
\langle z, w\rangle \geq s\|w\|^2 \quad \text{for all }\; z \in \partial^2_*\varphi(x,v)(w),\; x \in O, \;v\in\partial\varphi(x),\;w\in\R^n.
\end{equation}
The converse implication also holds under any one of the following conditions:
\begin{itemize}
    \item[\bf (i)] $\varphi$ is locally Lipschitz on $O$, and 
    \begin{equation}\label{ass:add}
\partial_p \varphi(x) = \partial\varphi(x) \quad \text{for all }\; x \in O.
\end{equation}
\item[\bf (ii)] $\varphi$ is hypoconvex on $O$, and $O = \R^n$. 
\end{itemize}
\end{Theorem}
\begin{proof}  It follows from \cite[Theorem~3.6]{roc24} that for each $x \in O, \; v \in \partial\varphi(x)$, the modulus of variational convexity can be expressed in terms of \eqref{eq:2ndstrict} as follows:
\begin{equation}\label{eq:2ndviavariational-2}
\operatorname{cnv}\ph(x|v)  = \operatorname{def} \partial^2_*\varphi(x,v),
\end{equation}
where 
$$
\operatorname{def} \partial^2_*\varphi(x,v) := \inf \{\langle z,u\rangle |\;  z \in \partial^2_*\varphi(x,v)(u),\; \|u\|=1\},
$$
under the assumption that $v\in\partial_p \varphi(x)$.  
We first verify the implication ``($\Longrightarrow$)''. Indeed, suppose that $\varphi$ is $s$-convex on $O$, then $\varphi$  is variationally $s$-convex at any point $x \in O$, which implies that $\operatorname{cnv}\varphi(x|v) \geq s$ for all $x \in O$ and  $v\in \partial\varphi(x)$. Pick any $x \in O$, $v \in \partial\varphi(x)$. Note that $v \in \partial_p \varphi(x)$ due to the variational $s$-convexity at $x$ for $v$. Combining the latter with  \eqref{eq:2ndviavariational-2}, we obtain 
\begin{equation}\label{2ndnorm1-2}
\langle z, u\rangle \geq s  \quad \text{for all }\; z \in (D_*\partial\varphi)(x,v)(u), \; \|u\|=1. 
\end{equation}
Pick any $w \in \R^n \setminus \{0\}$, $z \in (D_*\partial\varphi)(x,v)(w)$, and define $u:= w/\|w\|$. It follows from the fact that $(D_*\partial\varphi)(x,v)$ is positively homogeneous  that $z/\|w\| \in (D_*\partial\varphi)(x,v)(u)$.   Combining the latter with the fact that $\|u\|=1$ and \eqref{2ndnorm1-2}, we must have 
$$
\left\langle \frac{z}{\|w\|}, u  \right\rangle \geq s \iff \left\langle \frac{z}{\|w\|}, \frac{w}{\|w\|}  \right\rangle \geq s \iff \langle z , w\rangle \geq s\|w\|^2,
$$
which verifies \eqref{eq:graphbased} in the case where $w \ne 0$. In the case where $w =0$, \eqref{eq:graphbased} holds trivially. 

\medskip\noindent
($\Longleftarrow$):  We now assume that \eqref{eq:graphbased} holds and that one of conditions {\bf (i)} or {\bf (ii)} is satisfied. We aim to show that $\varphi$ is $s$-convex on $O$.  
For every $x \in O$ and $v \in \partial\varphi(x)$, we also have $v \in \partial_p \varphi(x)$ when {\bf (i)} or {\bf (ii)} holds. Therefore, 
it follows immediately from \eqref{eq:2ndviavariational-2} that $\operatorname{cnv}\varphi(x|v) \geq s$. Using Corollary~\ref{co:mainre} and Theorem~\ref{thm:hypo}, we deduce that $\varphi$ is $s$-convex on $O$. The proof is complete. 
\end{proof}
Note that Condition~\eqref{ass:add} can be implied by the prox-regularity of $\varphi$ on $O$, but the converse does not necessarily hold. This can be shown in the following example.  
\begin{Example}\rm 
Consider the function  $\varphi:\R\to\R$ by
$$
\varphi(x)=
\begin{cases}
|x|+|x|^3\sin(|x|^{-3/2})&  \text{if } x\neq 0\\
0,&\text{otherwise}.
\end{cases}
$$
Then $\varphi$ is locally Lipschitz on $\R$ and satisfies
$$
\partial_p\varphi(x)=\partial\varphi(x)
\qquad \text{for all }x\in\R,
$$
but $\varphi$ is not prox-regular at $0$ for $v=1$. 
Indeed, for every $x\neq 0$, the function $\varphi$ is smooth in a neighborhood of $x$. Hence,
$$
\partial_p\varphi(x)=\partial\varphi(x)
=\{\varphi'(x)\}.
$$
Moreover, we have  $
\partial_p\varphi(0)
=\partial\varphi(0)
=[-1,1]$. 
Therefore, the proximal and limiting subdifferentials coincide at every point of $\R$. 
We next show that $\varphi$ is not prox-regular at $0$ for $v=1$. For $x>0$, direct computation gives
$$
\varphi'(x)
=
1+3x^2\sin(x^{-3/2})
-\frac32x^{1/2}\cos(x^{-3/2}),
$$
and
$$
\varphi''(x)
=
6x\sin(x^{-3/2})
-\frac{21}{4}x^{-1/2}\cos(x^{-3/2})
-\frac94x^{-2}\sin(x^{-3/2}).
$$
Consider the sequence
$$
x_k=\left(\frac{\pi}{2}+2k\pi\right)^{-2/3}
\quad \text{for all }\; k\in\N.
$$
Then
$$
x_k\downarrow 0\quad \text{and }\; 
v_k:=\varphi'(x_k)=1+3x_k^2\to 1,
$$
and
$$
\varphi''(x_k)
=6x_k-\frac94x_k^{-2}
\to -\infty \; \text{as }\; k \to \infty.
$$
For an arbitrary $w\in\R$, choose $
w_k=w+\sqrt{x_k}$ for $k \in \N$. 
Then $w_k\to w$ and $w_k\neq 0$ for all sufficiently large $k$. By the second-order Taylor expansion at $x_k$, we can choose $t_k>0$ sufficiently small, with $t_k\downarrow0$, such that
$$
\frac{2\big[\varphi(x_k+t_kw_k)-\varphi(x_k)-t_kv_kw_k\big]}
{t_k^2}
\leq \varphi''(x_k)w_k^2+1.
$$
If $w\neq0$, then $w_k^2\to w^2>0$, and hence $\varphi''(x_k)w_k^2\to -\infty$. If $w=0$, then $w_k^2=x_k$, and
$$
\varphi''(x_k)w_k^2
=6x_k^2-\frac94x_k^{-1}
\to-\infty  \; \text{as }\; k \to \infty.
$$
Thus, in either case, the above quotient tends to $-\infty$. Consequently, $
d_*^2\varphi(0|1)(w)=-\infty$ for all $w \in \R$. 
Therefore, $\varphi$ is not prox-regular at $0$ for $v=1$. 
\end{Example}
We next provide another characterization of convexity in terms of quadratic bundles. This notion was originally introduced in \cite{roc} for convex functions and subsequently extended in \cite{gtdquad,quadcharvar,roc24} to a broader class of extended-real-valued functions $\varphi\colon\R^n\to\overline{\R}$. To proceed, we recall some relevant definitions.

\medskip 
A function $\varphi$ is said to be {\em twice epi-differentiable} at $\bar{x}$ for $\bar{v}$ if, for every $w\in\R^n$ and every sequence $\tau_k\downarrow0$, there exists a sequence $w^k\to w$ such that
$$
\Delta^2_{\tau_k}\varphi(\bar{x}|\bar{v})(w^k)
\to d^2\varphi(\bar{x}|\bar{v})(w)
\quad\text{as }k\to\infty.
$$
Given a symmetric matrix $A\in\R^{n\times n}$, we denote by
$$
q_A(w):=\langle w,Aw\rangle 
\quad \text{for all }\; w\in\R^n,
$$
the quadratic form associated with $A$. A function $q\colon\R^n\to\overline{\R}$ is called a {\em generalized quadratic form} if it can be represented as $
q=q_A+\delta_L$, 
where $L$ is a linear subspace of $\R^n$ and $A\in\R^{n\times n}$ is symmetric.
We say that $\varphi\colon\R^n\to\overline{\R}$ is {\em generalized twice differentiable} at $\bar{x}$ for $\bar{v}\in\partial\varphi(\bar{x})$ if it is twice epi-differentiable at $\bar{x}$ for $\bar{v}$ and its second-order subderivative $d^2\varphi(\bar{x}|\bar{v})$ is a generalized quadratic form. 
The {\em quadratic bundle} of $\varphi$ at $\bar{x}\in\operatorname{dom}\varphi$ for $\bar{v}\in\partial\varphi(\bar{x})$, denoted by $\operatorname{quad}\varphi(\bar{x}|\bar{v})$, is defined as the collection of all generalized quadratic forms $q$ for which there exists a sequence $
(x_k,v_k)\xrightarrow{\Omega_\varphi}(\bar{x},\bar{v}), 
\varphi(x_k)\to\varphi(\bar{x})$ 
such that the sequence of generalized quadratic forms $
q_k:=\frac{1}{2}d^2\varphi(x_k| v_k)$ 
converges epigraphically to $q$ as $k\to\infty$, where
$$
\Omega_\varphi:=
\left\{
(x,v)\in\operatorname{gph}\partial\varphi
\;\middle|\;
\varphi\text{ is generalized twice differentiable at }x
\text{ for }v
\right\}.
$$
Now we are ready to obtain the following characterization of convex functions. 
\begin{Theorem}[characterization   via  quadratic bundles]\label{thm:quad} Let $\varphi:\R^n\to\overline{\R}$ be a proper, l.s.c. function, $O \subset \R^n$ be an open, convex set, and let $s \in \R$.   If $\varphi$ is $s$-convex on $O$, then the following condition holds: 
\begin{equation}\label{eq:quadbased} 
q(w) \geq \frac{s}{2} \|w\|^2 \quad \text{for all }\; q \in \operatorname{quad}\varphi(x|v),\; x \in O, \;v\in\partial\varphi(x),\;w\in\R^n.
\end{equation}
The converse implication also holds under any one of the following conditions:
\begin{itemize}
    \item[\bf (i)]  $\varphi:\R^n \to \oR$ be locally Lipschitz   prox-regular on $O$. 
    \item[\bf (ii)] $\varphi$ is hypoconvex on $O$, and $O=\R^n$.  
\end{itemize}
\end{Theorem}
\begin{proof}
We first verify the implication ``($\Longrightarrow$)''. Indeed, suppose that $\varphi$ is $s$-convex on $O$. Fix any $x \in O$ and $v \in \partial\varphi(x)$, then $\varphi$  is variationally $s$-convex at $x$ for $v$. We use \cite[Theorem~4.4]{quadcharvar} to deduce that  \eqref{eq:quadbased} holds.  To verify the reverse implication    ``($\Longleftarrow$)'', we suppose that \eqref{eq:quadbased} holds, and that one of conditions {\bf (i)} or {\bf (ii)} is satisfied. Fix any $\mu < s$, $x \in O$, and $v \in \partial\varphi(x)$,  by using again  \cite[Theorem~4.4]{quadcharvar}, it follows that $\varphi$ is  variationally $\mu$-convex at $x$ for $v$. Therefore, $\operatorname{cnv}\varphi(x|v) \geq \mu$.  Due to Corollary~\ref{co:mainre} and Theorem~\ref{thm:hypo}, $\varphi$ is $\mu$-convex on $O$, which means that 
$$
\varphi(tx + (1-t)y) \leq t \varphi(x) + (1-t)\varphi(y) - \frac{\mu}{2}t(1-t)\|x-y\|^2 \quad \text{for all }\; x, y\in O, t \in [0,1].
$$
Since $\mu<s$ is chosen arbitrarily, we can let $\mu \to s^-$ to obtain 
$$
\varphi(tx + (1-t)y) \leq t \varphi(x) + (1-t)\varphi(y) - \frac{s}{2}\|x-y\|^2  \quad \text{for all }\; x, y\in O, t \in [0,1],
$$
and thus $\varphi$ is $s$-convex on $O$. The proof is complete. 
\end{proof} 

\section{Convexity of Sets via Variational Modulus of Convexity}\label{sec:localset}
In this section, we study the relationship between the classical convexity of a closed set and the variational convexity of its indicator function. We then establish a characterization of convex sets in terms of strict second subderivatives. We begin by presenting a result that connects the convexity of a closed set with the variational convexity of its indicator function.

\begin{Theorem}[local convexity of sets and variational convexity of indicator functions]\label{thm:set}
Let $\Omega\subset\R^n$ be a nonempty closed set and
let $\bar{x}\in\Omega$. Then $\Omega$ is locally convex around $\bar{x}$ if and only if the indicator function
$\delta_\Omega$ is variationally convex at $\bar{x}$ for
$\bar{v}:=0$.
\end{Theorem} 
\begin{proof} If $\Omega$ is locally convex around $\ox$, then its indicator
function $\delta_\Omega$ is locally convex around $\ox$.
Consequently, $\delta_\Omega$ is variationally convex at $\ox$
for $\ov=0$.

\medskip We next  verify the reverse implication. Suppose that $\delta_\Omega$ is variationally convex at $\ox$
for $\ov=0$.  
Since $\delta_\Omega$ is variationally convex at $\bar{x}$
for $\bar{v}=0$, by the definition of variational convexity,
there exist open convex neighborhoods $U$ of $\bar{x}$ and
$V$ of $0$, a number $\varepsilon>0$, and a proper l.s.c. convex function
$g:U\to\overline{\R}$ such that $
g(x)\leq \delta_\Omega(x)$ for all  $x\in U$ 
and
\begin{equation}\label{eq:VC-graph}
\big(U_\varepsilon\times V\big)
\cap\operatorname{gph}\partial \delta_\Omega
=
\big(U\times V\big)
\cap\operatorname{gph}\partial g, \quad g(x) = \delta_\Omega(x) \text{ at common elements}\;(x,v), 
\end{equation}
where
$$
U_\varepsilon
:=
\left\{x\in U:
\delta_\Omega(x)<\delta_\Omega(\bar{x})+\varepsilon
\right\} = \left\{x\in U:
\delta_\Omega(x)<\varepsilon
\right\}.
$$
By the definitions of $U_\epsilon$ and $\delta_\Omega$, it follows that $
U_\varepsilon=\Omega\cap U$.  
We claim that
\begin{equation}\label{eq:Omega-min}
\Omega\cap U
=
\operatorname*{argmin}_{u\in U}g(u).
\end{equation}
To prove the inclusion ``$\subset$'', take any
$x\in\Omega\cap U$. Since $0\in N_\Omega(x) = \partial\delta_\Omega(x)$, we have
$$
(x,0)\in
(U_\varepsilon\times V)\cap
\operatorname{gph}\partial\delta_\Omega.
$$
Combining the latter with \eqref{eq:VC-graph}, we deduce that $
(x,0)\in\operatorname{gph}\partial g$, or $0\in\partial g(x)$. 
Since $g$ is convex, the subgradient inequality gives
$$
g(y)\geq g(x)+\langle 0,y-x\rangle
=g(x)
\qquad\text{for all }y\in U,
$$
which means that $\displaystyle x \in \operatorname*{argmin}_{u\in U}g(u)$, and thus $\displaystyle \Omega\cap U
\subset 
\operatorname*{argmin}_{u\in U}g(u).$ 
For the reverse inclusion ``$\supset$'', let
$\displaystyle x\in\operatorname*{argmin}_{u\in U}g(u)$. It follows that $g(y) \geq g(x)$ for all $y \in U$, which implies that
$$
\liminf_{y \to x} \frac{g(y)-g(x)-\langle 0, y-x\rangle}{\|y-x\|} \geq 0,
$$
and thus $0 \in \widehat{\partial}g(x)\subset \partial g(x)$. 
Combining the latter with the fact that $0\in V$, we obtain 
$$
(x,0)\in
(U\times V)\cap\operatorname{gph}\partial g.
$$
By \eqref{eq:VC-graph},
we obtain $
(x,0)\in
(U_\varepsilon\times V)\cap
\operatorname{gph}\partial\delta_\Omega.$ 
It follows that $
x\in U_\varepsilon=\Omega\cap U,$ and thus  $\displaystyle \Omega\cap U
\supset 
\operatorname*{argmin}_{u\in U}g(x)$.
Therefore, \eqref{eq:Omega-min} is verified. 
Finally, since $g$ is convex on the convex set $U$,
its minimizer set over $U$ is convex. 
Combining this with \eqref{eq:Omega-min}, we deduce that  $\Omega\cap U$ is convex.
This proves that $\Omega$ is locally convex at
$\bar{x}$.
\end{proof}
Unlike the equivalence  between the convexity and variational convexity of
functions established in Theorem~\ref{tm:equivalent}, Theorem~\ref{thm:set}
requires only the variational convexity of the indicator function at
$\ox$ for the zero subgradient $\ov=0$. In particular, variational
convexity at $\ox$ for nonzero subgradients is not needed to guarantee
the local convexity of the set. This distinction is noteworthy when
compared with the corresponding result for functions. The following
example illustrates that variational convexity of the indicator
function at $\ox$ for a nonzero subgradient, without variational
convexity for the zero subgradient, is insufficient to guarantee
local convexity around $\ox$.
\begin{Example}\rm Consider the unit circle  in $\R^2$ 
$$
\Omega:= \{x\in \R^2|\; \|x\|=1\}= \{(x_1,x_2) \in \R^2|\; x_1^2+x_2^2 =1\},
$$
and let $\ox:=(1,0)$. We show that $\delta_\Omega$ is
variationally convex at $\ox$ for some nonzero subgradient,
although $\Omega$ is not locally convex around $\ox$.  Indeed, we have 
$$
\partial\delta_\Omega(x)= N_\Omega(x)=\{\lambda x\mid \lambda\in\R\}
\quad \text{for all }\;  x\in\Omega.
$$
In particular, for any $a>0$, setting $
\ov:=a\ox=(a,0)$, 
we have $\ov\in\partial\delta_\Omega(\ox)$.
We next verify that $\delta_\Omega$ is variationally convex
at $\ox$ for $\ov$. Since $\|x\|=1$ for every $x\in\Omega$,
whenever $
x\in\Omega,  v\in\partial\delta_\Omega(x),$
there exists $\lambda\in\R$ such that $v=\lambda x$.
Moreover, $
\lambda=\langle v,x\rangle.$ 
Since $\langle\ov,\ox\rangle=a>0$, the continuity of the
inner product guarantees the existence of  neighborhoods
$U$ of $\ox$ and $V$ of $\ov$ such that
$$
\langle v,x\rangle>\frac{a}{2}
\qquad
\text{for all }x\in U\cap\Omega,\ v\in V.
$$
Consequently, for every $
(x,v)\in
(U\times V)
\cap\operatorname{gph}\partial\delta_\Omega$, 
we can write $v=\lambda x$, where
$\lambda=\langle v,x\rangle>\frac{a}{2}.$ 
Now, let $x'\in\Omega$. Using $\|x\|=\|x'\|=1$, we obtain
$$
\|x'-x\|^2
=\|x'\|^2+\|x\|^2-2\langle x,x'\rangle
=2-2\langle x,x'\rangle.
$$
Therefore,
$$
\langle x,x'\rangle-1
=-\frac{1}{2}\|x'-x\|^2.
$$
It follows that
$$
\langle v,x'-x\rangle
=\lambda\langle x,x'-x\rangle
=\lambda(\langle x,x'\rangle-1)
=-\frac{\lambda}{2}\|x'-x\|^2
\leq-\frac{a}{4}\|x'-x\|^2.
$$
Hence,
\[
\delta_\Omega(x')
\geq
\delta_\Omega(x)+\langle v,x'-x\rangle
+\frac{a}{4}\|x'-x\|^2.
\]
The inequality is automatically satisfied when $x'\notin\Omega$,
since then $\delta_\Omega(x')=+\infty$.
It follows from Proposition~\ref{subgradientvar} that  $\delta_\Omega$ is variationally convex at
$\ox$ for the nonzero subgradient $\ov$. 
Finally, we show that $\Omega$ is not locally convex around
$\ox$. Indeed, for   any sufficiently small $\theta>0$ such that $\cos^2 \theta <1$, consider
$$
x=(\cos\theta,\sin\theta),
\qquad
y=(\cos\theta,-\sin\theta).
$$
Both points belong to $\Omega$ and can be made arbitrarily
close to $\ox$. However,
$$
\frac{x+y}{2}=(\cos\theta,0)
\quad \text{and }\; 
\left\|\frac{x+y}{2}\right\|^2
=\cos^2\theta<1.
$$
Thus, $(x+y)/2\notin\Omega$, showing that $\Omega$
is not locally convex around $\ox$.
\end{Example}

Our next objective is to characterize the local convexity of a closed set in terms of the variational modulus of convexity and strict second subderivatives. To this end, we  recall some preliminary concepts. By Proposition~\ref{subgradientvar}, a proper, l.s.c. function $\varphi:\R^n\to \oR$ is variationally $s$-convex at $\ox$ for $\ov \in\partial\varphi(\ox)$ if and only if there exists $\epsilon>0$ such that 
\begin{equation}\label{eq:proxineq}
\varphi(x')\geq \varphi(x)+\langle v,x'-x\rangle
+\frac{s}{2}\|x'-x\|^2
\end{equation}
for all $x'  \in B(\ox,\epsilon)$, $(x,v) \in \gph T_\epsilon^\varphi$, where $T_\epsilon^\varphi$ is the $\varphi$-attentive $\epsilon$-localization of $\partial\varphi$ around $(\ox,\ov)$ defined by 
\begin{equation}\label{eq:localization}
\gph T_\epsilon^\varphi = \{(x,v) \in \gph\partial\varphi|\; (x,v) \in B(\ox,\epsilon)\times B(\ov,\epsilon), \; \varphi(x)<\varphi(\ox) +\epsilon \}. 
\end{equation}
Following \cite[Definition~3.3]{gtdquad}, we recall the notion of
the {\em radius of variational $s$-convexity}. More precisely,
$\varphi$ is said to be {\em variationally $s$-convex} at $\ox$
for $\ov$ with radius $\epsilon>0$ if \eqref{eq:proxineq} holds.

 The following lemma extends the formula established in
\cite[Lemma~3.3]{ar26}. While the result in \cite{ar26} provides a formula for the strict
second subderivative of the Moreau envelope only at a single
reference point, our result establishes this formula throughout
the localization neighborhood. 

\begin{Lemma}[strict second subderivatives of Moreau envelopes]\label{lem:generalized-second-subderivative}
Let $\varphi:\R^n\to \oR$ be proper, l.s.c., and  prox-bounded. Suppose that $\varphi$ is $r$-level prox-regular at $\bar{x}$ for
$\bar{v}\in\partial \varphi(\bar{x})$ with radius $\epsilon>0$.  Then \eqref{eq:proxineq} holds with $s = -r$. 
Moreover, for all sufficiently small $\lm\in (0,1/r)$, there exists   an open convex neighborhood $U_\lm$ of $\ox +\lm \ov$ such that $e_\lm \varphi$ is $\mathcal{C}^{1,1}$ on $U_\lm$ and 
\begin{equation}\label{gradientC1}
\nabla e_\lm\varphi (u) = \frac{1}{\lm}\left(u - P_\lm\varphi(u) \right) = (\lm I + (T_\epsilon^\varphi)^{-1})^{-1}(u) \quad \text{for all }\; u \in U_\lm,
\end{equation}
Then, for every $u \in U_\lm$, $\xi\in\R^n$, we have
\begin{equation}\label{eq:generalized-second-subderivative}
d_*^2(e_\lambda \varphi)(u|v)(\xi)
=e_{\lambda/2}\big[d_*^2\varphi(x| v)\big](\xi),
\end{equation}
where 
\begin{equation}\label{eq:xv}
x:=P_\lambda \varphi(u),
\qquad
v:=\nabla e_\lambda \varphi(u)
=\frac{u-x}{\lambda}.
\end{equation}
\end{Lemma}
\begin{proof} It follows from \cite[Lemma~5.2]{gtdquad} that for all sufficiently small $\lm\in(0,1/r)$, there exists   an open convex neighborhood $U_\lm$ of $\ox+\lm\ov$ such that $e_\lm\varphi$ is $\mathcal{C}^{1,1}$ on $U_\lm$ and \eqref{gradientC1} is satisfied. Fix such $\lm$, and pick any $u \in U_\lm$, $\xi \in \R^n$. We define $x, v$ as in \eqref{eq:xv}, and prove the equality in \eqref{eq:generalized-second-subderivative}
by establishing both inequalities.  

\medskip We first prove ``$\geq$'' in \eqref{eq:generalized-second-subderivative}. 
Consider an arbitrary sequence $
(u_k,w_k,\tau_k,\xi_k)
\to (u,v,0,\xi)$ as $k \to \infty$, $w_k \in \partial e_\lm \varphi(u_k)$, $\tau_k \geq 0$ for all $k \in \N$.  
Since $e_\lambda \varphi$ is continuously differentiable on
$U_\lambda$, we may assume, for all sufficiently large $k$, that $
w_k=\nabla e_\lambda \varphi(u_k)$. 
For such $k$ sufficiently large, we  define 
$$
x_k:=P_\lambda \varphi(u_k),
\qquad
v_k:=\nabla e_\lambda \varphi(u_k)=w_k.
$$
By \eqref{gradientC1}, we have $u_k=x_k+\lambda v_k$ and  $
(x_k,v_k)\in
\operatorname{gph} T_\epsilon^\varphi$ 
for all sufficiently large $k$. 
Furthermore,
$$
x_k= P_\lm \varphi(u_k)\to P_\lm \varphi(u)=x 
\quad\text{and }\; 
v_k = \nabla e_\lm \varphi(u_k)  \to \nabla e_\lm\varphi(u)= v.
$$
Recall that
$$
e_\lambda \varphi(u_k)
=\varphi(x_k)+\frac{1}{2\lambda}\|u_k-x_k\|^2
=\varphi(x_k)+\frac{\lambda}{2}\|v_k\|^2.
$$
Since $e_\lambda \varphi(u_k)\to e_\lambda \varphi(u)$ and $v_k\to v$,
we obtain $
\varphi(x_k)\to \varphi(x)$ as $k \to \infty$.
Now define
$$
\eta_k:=
\frac{P_\lambda \varphi(u_k+\tau_k\xi_k)-x_k}{\tau_k}.
$$
Since $P_\lambda \varphi$ is locally Lipschitz on $U_\lambda$,
the sequence $\{\eta_k\}$ is bounded. Passing to a subsequence,
we may assume that $
\eta_k\to\eta$
for some $\eta\in\R^n$. Applying the finite-scale identity (10) in   the proof of \cite[Lemma 3.3]{ar26} at
$(x_k,v_k)$, with $u_k=x_k+\lambda v_k$, we obtain 
\begin{equation}\label{eq:finite-scale-identity}
\Delta_{\tau_k}^2[e_\lambda\varphi](u_k|v_k)(\xi_k)
=
\inf_{\theta\in\R^n}
\left\{
\Delta_{\tau_k}^2\varphi(x_k|v_k)(\theta)
+\frac{1}{\lambda}\|\xi_k-\theta\|^2
\right\},
\end{equation}
We next show that the infimum in \eqref{eq:finite-scale-identity} is attained at $\theta = \eta_k$. 
Indeed, for each $\theta \in \R^n$, we set $z=x_k+\tau_k\theta$. We have the following equalities 
\begin{align*}
&\frac{\tau_k^2}{2}
\left[
\Delta_{\tau_k}^2\varphi(x_k| v_k)(\theta)
+\frac{1}{\lambda}\|\xi_k-\theta\|^2
\right]\\
&\quad=
\varphi(z)-\varphi(x_k)-\langle v_k,z-x_k\rangle
+\frac{1}{2\lambda}\|x_k+\tau_k\xi_k-z\|^2\\
&\quad=
\varphi(z)+\frac{1}{2\lambda}
\|z-(u_k +\tau_k\xi_k)\|^2
-\varphi(x_k)-\frac{\lambda}{2}\|v_k\|^2
-\tau_k\langle v_k,\xi_k\rangle.
\end{align*}
Since the last three terms are independent of $z$, minimizing
the expression over $\theta\in\R^n$ is equivalent to
minimizing
$$
z \mapsto \varphi(z)+\frac{1}{2\lambda}
\|z-(u_k +\tau_k\xi_k)\|^2
$$
over $z\in\R^n$. By the definition of the proximal
mapping, the latter minimization is attained at $
z=P_\lambda\varphi(u_k+\tau_k\xi_k)$.
Consequently, the infimum in \eqref{eq:finite-scale-identity} is attained at
$$
\theta=
\frac{P_\lambda\varphi(u_k+\tau_k\xi_k)-x_k}{\tau_k}=\eta_k.
$$ 
Since $\varphi(x_k)\to\varphi(x)$, $v_k\in\partial\varphi(x_k)$, and $\tau_k \geq 0$ for all sufficiently large $k$, then the sequence $
(x_k,v_k,\tau_k,\eta_k)
\to (x,v,0,\eta)$ 
is admissible in the definition of
$d_*^2\varphi(x| v)(\eta)$. Consequently, we have 
$$
\liminf_{k\to\infty}
\Delta_{\tau_k}^2\varphi(x_k| v_k)(\eta_k)
\geq d_*^2\varphi(x|v)(\eta).
$$
Combining the latter with \eqref{eq:finite-scale-identity} and the fact that the infimum is attained at $\eta_k$, we deduce that 
$$
\liminf_{k\to\infty}
\Delta_{\tau_k}^2[e_\lambda \varphi](u_k|v_k)(\xi_k)
 \geq
d_*^2\varphi(x|v)(\eta)
+\frac{1}{\lambda}\|\xi-\eta\|^2 
  \geq
\alpha,
$$
where $$
\alpha:=
\inf_{\zeta\in\R^n}
\left\{
d_*^2\varphi(x|v)(\zeta)+\frac{1}{\lambda}\|\xi-\zeta\|^2
\right\}.
$$
Since the original admissible sequence was arbitrary, we conclude
that
\begin{equation}\label{eq:generalized-lower}
d_*^2(e_\lambda \varphi)(u|v)(\xi)
\geq\alpha,
\end{equation}
which verifies ``$\geq$'' in \eqref{eq:generalized-second-subderivative}. 

\medskip We now prove ``$\leq$'' in \eqref{eq:generalized-second-subderivative}. It follows from \eqref{eq:proxineq} and the definition of the strict second subderivative \eqref{strict2nd} that
\begin{equation}\label{eq:strict-second-lower-bound}
d_*^2\varphi(x|v)(\zeta)\geq-r\|\zeta\|^2
\qquad\text{for all }\zeta\in\R^n.
\end{equation}
Since $1/\lambda>r$, the function $
\zeta\mapsto
d_*^2\varphi(x| v)(\zeta)
+\frac{1}{\lambda}\|\xi-\zeta\|^2$
is coercive and bounded below. In particular, $\alpha$ is finite. 
Let $\delta>0$ be arbitrary. By the definition of $\alpha$, choose
$\bar{\zeta}\in\R^n$ such that
\begin{equation}\label{eq:approx-infimum}
d_*^2\varphi(x|v)(\bar{\zeta})
+\frac{1}{\lambda}\|\xi-\bar{\zeta}\|^2
\leq\alpha+\delta.
\end{equation}
By the definition in \eqref{strict2nd}, there exists
an admissible sequence $
(x_k,v_k,\tau_k,\zeta_k)
\to (x,v,0,\bar{\zeta})$  
such that $
\varphi(x_k)\to \varphi(x)$ 
and thus 
\begin{equation}\label{eq:realizing-sequence}
\liminf_{k\to\infty}
\Delta_{\tau_k}^2\varphi(x_k| v_k)(\zeta_k)
\leq d_*^2\varphi(x|v)(\bar{\zeta})+\delta.
\end{equation}
For each $k\in \N$, we define $
u_k:=x_k+\lambda v_k$ and 
$\xi_k:=\zeta_k+\xi-\bar{\zeta}$. 
Then $
u_k\to x+\lambda v=u$ and $
\xi_k\to\xi$ as $k \to \infty$. 
For all sufficiently large $k$, we deduce from \eqref{gradientC1} that  
$$
P_\lambda \varphi(u_k)=x_k \quad \text{and }\;  
\nabla e_\lambda \varphi(u_k)=v_k. 
$$
Moreover, we have 
$$
e_\lambda \varphi(u_k)
=\varphi(x_k)+\frac{1}{2\lambda}\|u_k-x_k\|^2
=\varphi(x_k)+\frac{\lambda}{2}\|v_k\|^2
\to
\varphi(x)+\frac{\lambda}{2}\|v\|^2
=e_\lambda \varphi(u).
$$
Thus $
(u_k,v_k,\tau_k,\xi_k)
\to (u,v,0,\xi)$ 
is admissible in the definition of
$d_*^2(e_\lambda \varphi)(u| v)(\xi)$. Applying the finite-scale identity (10) in   the proof of \cite[Lemma 3.3]{ar26} at
$(x_k,v_k)$, with $u_k=x_k+\lambda v_k$  and evaluating  at $\eta_k$, we obtain
$$
\Delta_{\tau_k}^2[e_\lambda \varphi](u_k|v_k)(\xi_k)
 \leq
\Delta_{\tau_k}^2\varphi(x_k| v_k)(\zeta_k)
+\frac{1}{\lambda}\|\xi_k-\zeta_k\|^2.
$$ 
Taking the lower limit and using
\eqref{eq:realizing-sequence}, we obtain
$$
d_*^2(e_\lambda \varphi)(u|v)(\xi)
 \leq
d_*^2\varphi(x| v)(\bar{\eta})+\delta
+\frac{1}{\lambda}\|\xi-\bar{\zeta}\|^2 
 \leq\alpha+2\delta,
$$
where the last inequality follows from
\eqref{eq:approx-infimum}. Since $\delta>0$ is arbitrary, we obtain
\begin{equation}\label{eq:generalized-upper}
d_*^2(e_\lambda \varphi)(u| v)(\xi)\leq\alpha.
\end{equation}
Combining \eqref{eq:generalized-lower} and
\eqref{eq:generalized-upper},
the proof is complete. 
\end{proof}

To derive our next results, the following remark is also useful.
\begin{Remark}[prox-boundedness from prox-regularity]\label{rem:prox-bounded} \rm Let $\varphi:\R^n\to \oR$ be $r$-level prox-regular at $\ox \in \dom\varphi$ for $\ov \in \partial\varphi(\ox)$ with radius $\epsilon>0$.  We define the function $\psi:\R^n \to \oR$ as follows: 
\begin{equation}\label{proxboundedness}
\psi(x):= \varphi (x)+\delta_{\mathbb{B}(\ox,\epsilon)}(x)\;\text{ whenever }\;x\in\R^n.
\end{equation} 
Then $\psi$ is prox-bounded and $r$-level prox-regular at $\ox$ for $\ov \in \partial\psi(\ox)$ with radius $\epsilon$. Furthermore, we have
\begin{equation}\label{1stchain}
\partial\psi(x) = \partial\varphi(x) \quad \text{for all }\; x \in B(\ox,\epsilon), 
\end{equation}
\begin{equation}\label{2ndchain}
d^2_* \psi(x|y) (u) = d^2_*\varphi(x,y)(u)\;\mbox{ for all }\;(x,y) \in\gph\partial \psi \cap (B(\ox,\epsilon)\times \R^n)\;\mbox{ and }\;u \in \R^n. 
\end{equation} 
\end{Remark}
We next establish a characterization of variational convexity in terms of  strict second subderivatives within a localization neighborhood. This result will play an important role in deriving the main results of this section.
\begin{Proposition}[neighborhood   characterizations for variational convexity via strict second subderivatives]\label{prop:neigh} Let $\varphi:\R^n \to \oR$ be prox-regular at $\ox \in \dom\varphi$ and $\ov \in \partial\varphi(\ox)$ with radius $\epsilon>0$. Then $\varphi$ is variationally convex at $\ox$ for $\ov$ with radius $\epsilon>0$ if and only if 
\begin{equation}\label{eq:2ndneighbor}
d_*^2 \varphi(x|v)(w) \geq 0 \quad \text{for all }\; (x,v) \in \gph T_\epsilon^\varphi, \; w \in \R^n,
\end{equation}
where $T_\epsilon^\varphi$ is defined as in \eqref{eq:localization}. 
\end{Proposition}
\begin{proof}  Suppose that $\varphi$ is variationally convex at $\ox$ for $\ov$ with radius $\epsilon>0$, we verify \eqref{eq:2ndneighbor}.  Indeed, the imposed assumption implies that \eqref{eq:proxineq} holds with $s=0$. Fix $(x,v) \in \gph T_\epsilon^\varphi$, and  $w \in \R^n$. Consider an arbitrary sequence $
(x_k,v_k,\tau_k,w_k)
\to (x,v,0,w)$  
such that $\tau_k \geq 0$, $v_k \in \partial\varphi(x_k)$ for all $k \in \N$, and $
\varphi(x_k)\to \varphi(x)$ as $k \to \infty$.  Then we have $(x_k,v_k) \in \gph T_\epsilon^\varphi$ for all sufficiently large $k \in \N$, and we deduce from \eqref{eq:proxineq} with $s=0$ that
$$
\Delta_{\tau_k}^2 \varphi(x_k|v_k)(w_k)= \frac{\varphi(x_k+\tau_k w_k)-\varphi(x_k)-\tau_k\langle v_k, w_k\rangle}{\frac{1}{2}\tau_k^2} \geq 0
$$
for all sufficiently large $k \in \N$.  
Since the original admissible sequence was arbitrary, we conclude
that $
 d_s^2 \varphi(x|v)(w) \geq 0$, which verifies \eqref{eq:2ndneighbor}.  

\medskip We next suppose that \eqref{eq:2ndneighbor} holds, we show that $\varphi$ is variationally convex at $\ox$ for $\ov$ with radius $\epsilon>0$. 
By Remark~\ref{rem:prox-bounded}, the prox-regularity of $\varphi$ at $\ox$ for $\ov$ with radius $\epsilon>0$ implies that the function $\psi$ defined as in \eqref{proxboundedness} is prox-bounded and prox-regular at $\ox$ for $\ov$ with radius $\epsilon>0$, and \eqref{1stchain} and \eqref{2ndchain} hold. 
It follows from Lemma~\ref{lem:generalized-second-subderivative} that for all sufficiently small $\lm >0$, there exists   an open convex neighborhood $U_\lm$ of $\ox +\lm \ov$ such that $e_\lm \psi$ is $\mathcal{C}^{1,1}$ on $U_\lm$ and 
\begin{equation}\label{eq:2ndofpsi}
d^2_* e_\lm \psi(u|v)(w) = e_{\lm/2} [d^2_*\psi(x|v)](w) 
\end{equation}
for all $u \in U_\lm$, $w \in \R^n$, and $x, v$ defined by 
$$
x:=P_\lambda \psi(u),\; 
v:=\nabla e_\lambda \psi(u)
=\frac{u-x}{\lambda}\; \text{and } (x,v) \in \gph T_\epsilon^\psi,
$$
where 
$$
\gph T_\epsilon^\psi = \{(x,v) \in \gph\partial\psi|\; (x,v) \in B(\ox,\epsilon)\times B(\ov,\epsilon), \; \psi(x)<\psi(\ox) +\epsilon \}.
$$
Note that $\gph T_\epsilon^\psi = \gph T_\epsilon^\varphi$ due to \eqref{1stchain}. 
Combining the latter with \eqref{2ndchain},  \eqref{eq:2ndneighbor} and \eqref{eq:2ndofpsi}, we obtain 
$$
d^2_* e_\lm \psi (u|v)(w) \geq 0 \quad \text{for all }\; u \in U_\lm, \; v  = \nabla e_\lm \psi(u), \; w \in \R^n.
$$
It follows from Theorem~\ref{convexvia2nd} that $e_\lm\psi$ is convex on $U_\lm$. By an argument similar to that used in the proof of \cite[Theorem~3.2]{kmp22convex}, we can show  that $\psi$ is variationally convex at $\ox$ for $\ov$ with radius $\epsilon$. Consequently, the same property holds for $\varphi$. This completes the proof.
    
\end{proof}

We are now ready to establish a characterization of locally convex sets in terms of the variational modulus of convexity and strict second subderivatives. 
\begin{Theorem}[characterizations  via variational modulus of convexity and strict second subderivatives]\label{thm:setchar} Let $\Omega \subset \R^n$ be a nonempty closed set. Then the following assertions are equivalent:
\begin{itemize}
    \item[\bf (i)] $\Omega$ is a  locally convex set.
    \item[\bf (ii)] $\operatorname{cnv}\delta_\Omega(x|v) \geq 0$ for all $x \in \partial\Omega$ and $v \in N_\Omega(x)$. 
    \item[\bf (iii)] $d^2_*\delta_\Omega(x|v)(w)\geq 0 $  for all $x \in \partial\Omega$, $v \in N_\Omega(x)$ and $w\in \R^n$. 
\end{itemize}
If in addition that $\Omega$ is connected, then $\Omega$ is convex if one of the above assertions holds. 
\end{Theorem}
\begin{proof} Suppose that {\bf (i)} holds. Then $\delta_\Omega$ is locally convex at every point $x\in\Omega$, and hence variationally convex at every such point. Consequently, {\bf (ii)} holds. The implication [{\bf (ii)}$\Longrightarrow${\bf (iii)}] follows from \eqref{cnv2nd}. We only need to verify the  implication [{\bf (iii)}$\Longrightarrow${\bf (i)}]. Suppose that the assertion {\bf (iii)} is satisfied. Fix any $\ox \in \partial\Omega$, by {\bf (iii)} and the fact that $d^2_*\delta_\Omega(x|v)= 0$ for any $x \in \text{\rm int}\Omega$, $v \in N_\Omega(x) = \{0\}$,  we deduce that for any $\epsilon>0$, we have 
\begin{equation}\label{eq:indicatorneigh}
d^2_*\delta_\Omega(x|v)(w) \geq 0 \quad \text{whenever }\; (x,v) \in \gph N_\Omega \cap (B(\ox,\epsilon)\times B(0,\epsilon)) = \gph T_\epsilon^{\Omega},
\end{equation}
where 
$$
\gph T_\epsilon^{\Omega}:= \{(x,v) \in \gph N_\Omega|\;  (x,v) \in B(\ox,\epsilon)\times B(0,\epsilon), \delta_\Omega(x) < \delta_\Omega(\ox) +\epsilon\}
$$
Moreover, the condition $d^2_*\varphi(\ox|0)\geq 0$ also implies that $\delta_\Omega$ is prox-regular at $\ox$ for $0$ with some radius $\epsilon>0$. Combining the latter with \eqref{eq:indicatorneigh}, we deduce from Proposition~\ref{prop:neigh} that $\delta_\Omega$ is variationally convex at $\ox$ for $0$. Using Theorem~\ref{thm:set}, we conclude that $\Omega$ is locally convex around $\ox$. Since $\ox$ is an arbitrary boundary point, $\Omega$ is a locally convex set. Thus, the assertion {\bf (i)} is verified. Therefore, the assertions {\bf (i)}, {\bf (ii)}, and {\bf (iii)} are equivalent.   When $\Omega$ is connected, its local convexity implies its global convexity, as established by Tietze in \cite{T36}. The proof is complete.  
    
\end{proof}

Consequently, we obtain the following characterizations of convexity for closed sets in terms of the strict graphical derivative \eqref{strictgraphical}, second-order subdifferentials \eqref{2nd}, and quadratic bundles as follows.

\begin{Corollary}[\rm characterization via strict graphical derivatives]\label{cor:strict} Let $\Omega\subset\R^n$ be a nonempty closed set. Then all the assertions of Theorem~\ref{thm:setchar} imply the following condition:
\begin{itemize}
    \item[{\bf (iv)}] $\langle z,w\rangle\geq0$ for all
    $z\in(D_*N_\Omega)(x,v)(w)$, 
    $x\in\partial\Omega$, $v\in N_\Omega(x)$, and $w\in\R^n$.
\end{itemize}
Moreover, the converse implication also holds provided that
\begin{equation}\label{eq:regular}
N_\Omega^P(x)=N_\Omega(x)
\quad\text{for all }x\in\partial\Omega.
\end{equation}
Consequently, under \eqref{eq:regular}, the assertions of Theorem~\ref{thm:setchar} are equivalent to {\bf (iv)}.
\end{Corollary}
\begin{proof} Since assertions {\bf (i)}, {\bf (ii)}, and {\bf (iii)} of Theorem~\ref{thm:setchar} are equivalent, it suffices to establish that {\bf (ii)}$\Longrightarrow${\bf (iv)} and, under condition~\eqref{eq:regular}, the converse implication {\bf (iv)}$\Longrightarrow${\bf (ii)}. Suppose that {\bf (ii)} holds, then for any $x \in \partial\Omega$ and $v \in N_\Omega(x)$,  the function $\delta_\Omega$ is prox-regular at $x$ for $v$. Using   \cite[Theorem~3.6]{roc24}, for each $x \in \partial\Omega$ and $v \in N_\Omega(x)$, the modulus of variational convexity can be expressed in terms of \eqref{eq:2ndstrict} as follows:
\begin{equation}\label{eq:2ndviavariational-3}
\operatorname{cnv}\delta_\Omega(x|v)  = \operatorname{def} \partial^2_*\delta_\Omega(x,v) = \operatorname{def}(D_*N_\Omega)(x,v),
\end{equation} 
where the last equality comes from the subdifferential continuity of indicator function $\delta_\Omega$. Combining {\bf (ii)} and \eqref{eq:2ndviavariational-3}, we obtain {\bf (iv)}. Now, we suppose that {\bf (iv)} holds, and condition \eqref{eq:regular} is satisfied. This condition is equivalent to the fact that  $\partial_p\delta_\Omega(x)= \partial \delta_\Omega(x)$ for all $x \in \partial\Omega$. It follows from  \cite[Theorem~3.6]{roc24} that \eqref{eq:2ndviavariational-3} holds, and thus $\operatorname{cnv}\delta_\Omega(x|v)\geq 0$ for all $x \in \partial \Omega$, $v\in N_\Omega(x)$, which verifies {\bf (ii)}. The proof is complete. 

\end{proof}

\begin{Corollary}[characterizations via second-order subdifferentials and quadratic bundles] Let $\Omega\subset\R^n$ be a nonempty closed set. Then all the assertions of Theorem~\ref{thm:setchar} imply the following conditions:
\begin{itemize}
    \item[{\bf (v)}] $\langle z,w\rangle\geq0$ for all
    $z\in \partial^2\delta_\Omega(x,v)(w)$, 
    $x\in\partial\Omega$, $v\in N_\Omega(x)$, and $w\in\R^n$.
    \item[\bf (vi)] $q(w) \geq 0$ for all
    $q\in \operatorname{quad}\delta_\Omega(x|v)$, 
    $x\in\partial\Omega$, $v\in N_\Omega(x)$, and $w\in\R^n$.
\end{itemize}
Moreover, all assertions in Theorem~\ref{thm:setchar} are equivalent to {\bf (v)} and {\bf (vi)} if $\delta_\Omega$ is prox-regular at any $x \in \partial\Omega$. 
\end{Corollary}
\begin{proof} The arguments follow closely those of Corollary~\ref{cor:strict}, where the variational modulus of convexity is characterized in terms of second-order subdifferentials and quadratic bundles, using \cite[Theorem~3.8]{roc24} and \cite[Theorem~4.5]{roc24}, respectively. Therefore, we omit the proof of this corollary. 
    
\end{proof}

\begin{Remark}[on the characterizations of convex sets] \rm Characterizations of locally convex sets have been investigated in \cite{chieuyen} using the second-order subdifferential \eqref{2nd}. However, the results in \cite{chieuyen} rely on rather restrictive structural assumptions on the set $\Omega\subset\R^n$. For instance, they require $\Omega$ to have a $\mathcal{C}^2$-smooth boundary or to admit a local representation as the epigraph of a $\mathcal{C}^1$-smooth or piecewise $\mathcal{C}^2$ function. In contrast, our approach does not rely on any specific structural assumptions on $\Omega$. Instead, we introduce a new perspective by exploiting the variational convexity of the indicator function $\delta_\Omega$ and then using the second-order generalized derivatives to investigate. This approach offers a new viewpoint for studying local convexity of closed sets in variational analysis and optimization. 
\end{Remark}

\section{Conclusion}\label{sec:conclusion}
In this paper, we have investigated the relationship between variational convexity and classical convexity of functions and sets. By exploiting the variational modulus of convexity and generalized second-order derivatives, we have established new characterizations of classical convexity through the framework of variational convexity.

\medskip 
We have demonstrated that variational convexity provides an effective framework for characterizing classical convexity through the variational modulus of convexity under the local Lipschitz continuity or hypoconvexity of functions. This framework further yields derivative-based criteria for convexity under  mild assumptions. We have also extended this approach to characterize the local convexity of sets through the variational convexity of their indicator functions. In particular, our results apply to arbitrary closed sets without imposing additional geometric assumptions, such as boundary smoothness or specific local representations. These findings establish a novel connection between the geometric properties of sets and the second-order variational analysis of their indicator functions.

\section*{Declarations}

\section*{Statements and Declarations}

\begin{itemize}
\item \textbf{Funding:} Vo Thanh Phat was partially supported by the Early Career Scholars Program 2026 at the University of North Dakota under project number 43700-2375-UND0031286.

\item \textbf{Competing Interests:} The authors have no relevant financial or non-financial interests to disclose.

\item \textbf{Data Availability:} Data sharing is not applicable to this article as no datasets were generated or analysed during the current study.

\item \textbf{Author Contributions:} All authors contributed equally to this work. All authors read and approved the final manuscript.
\end{itemize}

\begin{thebibliography}{99}
\small 

\bibitem{ar26}  Adly, S.,     Rockafellar, R. T.: Variational modulus of convexity under Moreau envelopes: attenuation, recovery, and oracle tuning. HAL preprint hal-05563788, (2026), https://hal.science/hal-05563788 

\bibitem{ar27} Adly, S.,     Rockafellar, R. T.: Uniform quadratic growth and second-order variational curvature in Hilbert spaces.  HAL preprint hal-05768304, (2026), https://cnrs.hal.science/hal-05768304/

\bibitem{ChieuChuongYaoYen}   Chieu, N.H., Chuong, T.D., Yao, J.-C., Yen, N.D.: Characterizing convexity of a function by its Fr\'echet and limiting second-order subdifferentials. Set-Valued Var. Anal. 19, 75--96 (2011)

\bibitem{chieuyao11} Chieu, N.H.,  Yao, J.C.: Characterization of convexity for a piecewise $\mathcal{C}^2$ function
by the limiting second-order subdifferential. Taiwanese J. Math. 15, 31--42 (2011)

\bibitem{chieuhuy}  Chieu, N.H., Huy, N.Q.: Second-order subdifferentials and convexity of real-valued functions. Nonlinear Anal. 74, 154--160 (2011)

\bibitem{clmn16} Chieu, N.H.,  Lee, G.M.,  Mordukhovich, B.S.,  Nghia, T.T.A.: Coderivative characterizations of maximal monotonicity for set-valued mappings. J. Convex Anal. 23, 461--480 (2016)

\bibitem{clarke} Clarke, F.: Optimization and Nonsmooth Analysis. Wiley, New York (1983)


\bibitem{chieuyen} Chieu, N.H., Yen, N.D.: Convexity of sets and functions via second-order subdifferentials. Linear Nonlinear Anal. 5, 183--199 (2019)

\bibitem{dd22} Davis, D., Drusvyatskiy, D.: Proximal methods avoid active strict saddles of weakly convex functions. Found Comput Math. 22, 561--606 (2022). 

\bibitem{g25} Gfrerer, H. On second-order variational analysis of variational convexity of prox-regular functions. Set-Valued Var. Anal. 33, 8 (2025) 

\bibitem{ly10} Li, Y., Yeh, C.: Some characterizations of convex functions. Comput. Math. Appl. 59, 327--337 (2010)


\bibitem{kkmpbook} Khanh, P.D., Khoa, V. V. H., Mordukhovich, B.S., Phat, V.T.: Variational convexity of functions in Banach spaces, in Functional Analysis and Continuous Optimization, J. M. Amig\'o, M. J. C\'anovas et al. (eds.). Springer Proceedings in Mathematics and Statistics 424, pp.\ 237--260, Springer, Cham, Switzerland (2023)

\bibitem{KKMP23-2} Khanh, P.D., Khoa, V. V. H., Mordukhovich, B.S., Phat, V.T.:  Variational and strong variational convexity in infinite-dimensional variational analysis. SIAM J. Optim. 34, 2756--2787 (2024)

\bibitem{KKMP23-3} Khanh, P.D., Khoa, V. V. H., Mordukhovich, B.S., Phat, V.T.:  Second-order subdifferential optimality conditions in nonsmooth optimization. SIAM J. Optim. 35, 678--711 (2025) 

\bibitem{kkmpvn} Khanh, P.D., Khoa, V.V.H., Mordukhovich, B.S., Phat. V.T.: Local minimizers of nonconvex functions in Banach spaces via Moreau envelopes. Vietnam J. Math. 53, 803--813 (2025)

\bibitem{BorisKhanhPhat} Khanh, P.D., Mordukhovich, B.S., Phat, V.T.: A generalized Newton method for subgradient systems. Math. Oper. Res. 48, 1811--1845 (2023)


\bibitem{kmp22convex} Khanh, P.D., Mordukhovich, B.S., Phat, V.T.:  Variational convexity of functions and variational sufficiency in optimization. SIAM J. Optim. 33, 1121--1158 (2023) 

\bibitem{kmpmp} Khanh, P.D., Mordukhovich, B.S., Phat, V.T.:  Coderivative-based Newton methods in structured nonconvex and nonsmooth optimization. Math. Program. (2026). https://doi.org/10.1007/s10107-026-02420-7

\bibitem{kmptjogo} Khanh, P.D., Mordukhovich, B.S., Phat, V.T., Tran, D. B.:  Generalized damped Newton algorithms in nonsmooth optimization via second-order subdifferentials. J. Global Optim. 86, 93--122 (2023)

\bibitem{kmptmp} Khanh, P.D., Mordukhovich, B.S., Phat, V.T., Tran, D. B.:  Globally convergent coderivative-based generalized Newton methods in nonsmooth optimization. Math. Program. 205, 373--429 (2024)

\bibitem{kmptjogo25}  Khanh, P.D., Mordukhovich, B.S., Phat, V.T., Tran, D. B.:    Inexact proximal methods for
 weakly convex functions. J. Global Optim. 91, 611--646 (2025)

\bibitem{gtdquad}  Khanh, P.D., Mordukhovich, B.S., Phat, V.T., Viet, L. D.: Generalized twice differentiability and quadratic bundles in second-order variational analysis. Set-Valued Var. Anal 34, 11 (2026). https://doi.org/10.1007/s11228-026-00795-5

\bibitem{quadcharvar} Khanh, P.D., Mordukhovich, B.S., Phat, V.T.,  Viet, L. D.: Characterizations of variational convexity and tilt stability via quadratic bundles.  J. Convex Anal. 33, 903--924 (2026)

\bibitem{lwwy26} Luo, H., Wang, X., Wang, Z., Yang, X.: Level proximal subdifferential, variational convexity, and pointwise quadratic approximation. Math. Program. (2026). https://doi.org/10.1007/s10107-026-02355-z.


\bibitem{m92} Mordukhovich, B.S.: Sensitivity analysis in nonsmooth optimization. In: Field, D.A., Komkov, V.(eds) Theoretical Aspects of Industrial Design, pp.\ 32--46. SIAM Proc. Appl. Math. 58. Philadelphia, PA (1992)


\bibitem{Mordukhovich06} Mordukhovich, B.S.: Variational Analysis and Generalized Differentiation, I: Basic Theory, II: Applications. Springer, Berlin (2006)

\bibitem{Mor18} Mordukhovich, B.S.: Variational Analysis and Applications. Springer, Cham, Switzerland (2018)

\bibitem{Mor24} Mordukhovich, B.S.: Second-Order Variational Analysis in Optimization, Variational Stability, and Control: Theory, Algorithms, Applications. Springer, Cham, Switzerland (2024) 


\bibitem{phat26} Phat, V.T.: New characterizations of nonsmooth convex functions via generalized derivatives, 	arXiv:2603.24951 (2026)

\bibitem{poli90} Poliquin, R.A.: Subgradient monotonicity and convex functions. Nonlinear Anal. Theory Methods Appl. 14, 385--398 (1990) 


\bibitem{Rockafellar98} Rockafellar, R.T., Wets R.J-B.: Variational Analysis. Springer, Berlin (1998)

\bibitem{r76}  Rockafellar, R. T.: Monotone operators and the proximal point algorithm. SIAM J. Control Optim. 14, 877--898  (1976) 

\bibitem{r76mor}  Rockafellar, R. T.: Augmented Lagrangians and applications of the proximal point algorithm in convex programming. Math. Oper. Res. 1, 97--116 (1976)

\bibitem{r19} Rockafellar R.T.: Variational convexity and local monotonicity of subgradient mappings. Vietnam J. Math.  47, 547--561 (2019)

\bibitem{r19sv} Rockafellar, R.T.: Progressive decoupling of linkages in optimization and variational inequalities with elicitable convexity or monotonicity. Set-Valued Var. Anal. 27, 863--893 (2019)

\bibitem{r21} Rockafellar, R.T.: Characterizing firm nonexpansiveness of prox mappings both locally and globally. J. Nonlinear Convex Anal., 22(5) (2021) 

\bibitem{r21ppm} Rockafellar, R.T.: Advances in convergence and scope of the proximal point algorithm. J. Nonlinear Convex Anal. 22, 2347--2375 (2021)


\bibitem{roc} Rockafellar, R.T.: Augmented Lagrangians and hidden convexity in sufficient conditions for local optimality. Math. Program. 198, 159--194 (2023)

\bibitem{r23} Rockafellar, R.T.: Convergence of augmented Lagrangian methods in extensions beyond nonlinear programming. Math. Program. 199, 375--420 (2023)


\bibitem{varprox} Rockafellar, R.T.: 
Variational convexity and prox-regularity.  J. Convex Anal. 33, 1049--1060 (2026)

\bibitem{roc24} Rockafellar, R.T.:  Derivative tests for prox-regularity and the modulus of convexity.   J. Convex Anal. 33, 242--266 (2026)

 

\bibitem{r26} Rockafellar, R.T. Advances in progressive decoupling methodology for problem decomposition in optimization. Math. Program. (2026). https://doi.org/10.1007/s10107-026-02430-5



\bibitem{s81} Stromberg, K. R.: Introduction to Classical Real Analysis. Wadsworth International, Belmont, California (1981)

\bibitem{T36} Tietze, H.:  Über Konvexheit im kleinen und im grossen und über gewisse den Punkten einer Menge zugeordnete Dimensionszahlen. Math.~Z. 28, 697--707 (1938) 

\bibitem{wang10}  Wang, X.: On Chebyshev functions and Klee functions. J. Math. Anal. Appl.  368, 293--310 (2010)


\end{thebibliography}
\end{document}